\documentclass[12pt,twoside]{amsart}
\input{includeNice}
\usepackage{mabliautoref}
\usepackage{amssymb}
\usepackage{amscd}
\usepackage[abbrev,alphabetic]{amsrefs}
\usepackage{hyperref}
\usepackage[a4paper,margin=1in]{geometry}  
\usepackage{soul}

\usepackage{xypic}
\usepackage{tikz-cd}
\usepackage{caption}

\usepackage{url}
\RequirePackage{booktabs, multirow}
\RequirePackage{pgf}

\newcommand{\FB}[1]{\textcolor{blue}{(FB: #1)}}

\title[Rational points on 3-folds with nef anti-canonical class over $\mathbb{F}_q$]
{Rational points on 3-folds with nef anti-canonical class over finite fields} 
\author[F. Bernasconi and S. Filipazzi]{Fabio Bernasconi and Stefano Filipazzi} 
\subjclass[2020]{14E30, 14G05, 14G15, 14J20, 14J32}
\keywords{3-folds with nef anti-canonical class, rational points, finite fields}
\thanks{FB was partly supported by the grants $\#200021/169639$ and PZ00P2-21610 from the
Swiss National Science Foundation, SF was partly supported by ERC starting grant $\#804334$.}

\address{Dipartimento di Matematica “Guido
Castelnuovo”, SAPIENZA Università di Roma, Piazzale Aldo Moro 5, I-00185 
Roma}
\email{fabio.bernasconi@uniroma1.it}

\address{\'Ecole Polytechnique F\'ed\'erale de Lausanne, Chair of Algebraic Geometry  C3 625 
	(B\^atiment MA) Station 8 CH-1015 Lausanne 
} 
\email{stefano.filipazzi@epfl.ch}

\begin{document}
	
	\begin{abstract}
    We prove that a geometrically integral smooth projective 3-fold $X$ with nef anti-canonical class and negative Kodaira dimension over a finite field $\mathbb{F}_q$ of characteristic $p>5$ and cardinality $q=p^e > 19$ has a rational point. 
    Additionally, under the same assumptions on $p$ and $q$, we show that a smooth projective 3-fold $X$ with trivial canonical class and non-zero first Betti number $b_1(X) \neq 0$ has a rational point.
    Our techniques rely on the Minimal Model Program to establish several structure results for generalized log Calabi--Yau 3-fold pairs over perfect fields.
    \end{abstract}
	
	\maketitle
	\tableofcontents
	
\section{Introduction}

A fundamental principle in complex, algebraic, and arithmetic geometry states that the positivity properties of the canonical divisor $K_X$ of a smooth projective variety $X$ strongly influence the topology, geometry, and arithmetic of $X$. 
In this article, we investigate varieties with nef anti-canonical class, which can be regarded as being uniformly non-hyperbolic in some loose sense.

In the case $X$ is defined over the complex numbers $\mathbb{C}$, Cao and Horing \cite{CH19} proved that the universal cover of $X$ can be decomposed into a direct product consisting of copies of $\mathbb{C}$, along with simply connected $K$-trivial varieties and rationally connected ones.
From the point of view of Kobayashi hyperbolicity \cite{Kob98}, $X$ is then expected to be covered by holomorphic entire curves and its Kobayashi distance is expected to vanish everywhere \cite{Kob76}*{Problem F.2, p. 405}. 
On the arithmetic side, in characteristic 0, according to the predictions of the Green--Griffiths--Lang conjecture, rational points on varieties over number fields with $-K_X$ nef are expected to be potentially dense, i.e., the rational points are Zariski dense after passing to some finite field extension.

In this article, we aim to explore similar statements for varieties with $-K_X$ nef over finite fields and, more generally, over fields of prime characteristic $p>0$.
This has been motivated by a recent result of Ejiri and Patakfalvi \cite{EP23}, where the authors settle the Demailly--Peternell--Schneider conjecture in positive characteristic: They show that the Albanese morphism for these varieties is surjective.
Furthermore, under additional assumptions, they show that the Albanese morphism is locally isotrivial in the \'etale topology, similarly to what happens in characteristic 0 \cite{Cao19}.
In particular, this represents a first step towards a positive characteristic analog of the decomposition theorem in \cite{CH19}.
On the arithmetic side, when $k$ is a finite field, we propose the following generalization of a conjecture of Patakfalvi and Zdanowicz on the existence of rational points on $K$-trivial varieties over finite fields, cf. \cite{PZ}*{Conjecture 1.10}:

\begin{conjecture*} \label{conj: anti-nef-rat-pts}
    For every positive integer $d$ there exists an integer $q_d$ such that the following holds:
    If $X$ is a geometrically integral smooth projective variety with $-K_X$ nef  of dimension $d=\dim (X)$ over a finite field $\mathbb{F}_q$ and $q \geq q_d$, then $X(\mathbb{F}_q) \neq \emptyset$.
\end{conjecture*}

There are various results in the literature supporting \autoref{conj: anti-nef-rat-pts}.
When $X$ is a smooth Fano variety (or more generally a smooth rationally chain connected variety) over a finite field, Esnault proved that $X$ has a rational point \cite{Esn03}. 
If $X$ is geometrically an Abelian variety, it has a rational point by \cite{Lan55}.

The main result of this article is a proof of \autoref{conj: anti-nef-rat-pts} in the case of 3-folds of negative Kodaira dimension.

\begin{theorem}[Smooth case of \autoref{thm: intro_rat_pts_anti-nef}] \label{thm: intro-smooth}
    Let $p>5$ be a prime number and let $q \coloneqq p^e$ for some $e \geq 1$.
    Let $X$ be a geometrically integral smooth projective 3-fold over $\mathbb{F}_q$ with $-K_X$ is nef. 
    Suppose the Kodaira dimension of $X$ is $\kappa(X)= - \infty$.
    If $q>19$, then $X(\mathbb{F}_q) \neq \emptyset$ holds.
\end{theorem}

We are not able to fully address the existence of rational points in the case where $X$ is a smooth $K$-trivial 3-fold.
Nevertheless, using the recent results of \cite{EP23} together with a detailed study of the general fiber of the Albanese morphism, we can prove the following:

\begin{theorem}[Smooth case of \autoref{thm: intro_alb_non_trivial}]
 Let $p>5$ be a prime number and let $q \coloneqq p^e$ for some $e \geq 1$.
Let $X$ be a geometrically integral smooth projective 3-fold over $\mathbb{F}_q$ with $K_X \sim_{\mathbb{Q}} 0$.
If $b_1(X) \neq 0$ and $q>19$, then $X(\mathbb{F}_q) \neq \emptyset$ holds.
\end{theorem}

The restriction on the characteristic is due to the use of the Minimal Model Program (MMP) and vanishing theorems for 3-folds, which are known to be true for 3-folds in characteristic $p>5$ \cites{HX15, Bir16, BW17, BK23}.
As for the assumption on the cardinality $q>19$, this is the bound given for the existence of rational points on K3 surfaces over finite fields by the Weil conjectures (see \autoref{prop: rat_pts_reg_K_triv}).

\subsection{Singular case}

Our approach to the proof of \autoref{thm: intro-smooth} relies on several ideas and techniques of higher dimensional algebraic geometry.
In particular, we make extensive use of the MMP for 3-folds over perfect fields of characteristic $p > 5$, see \cites{GNT19}.
These techniques have already proved themselves fruitful in \cites{GNT19, NT20, Nak21}, where the MMP was used to establish Kodaira-type vanishing theorems for Witt vector cohomology of 3-fold with ample anti-canonical class.
In turn, these vanishing theorems allowed to generalize Esnault's result to the case of klt 3-fold Fano pairs.

We follow the same principle of systematically using the MMP to construct rational points on 3-folds with nef anti-canonical divisor. 
However, if $-K_X$ is only nef, we cannot reduce the problem of the existence of a rational point on $X$ to a purely cohomological computation using Witt vectors cohomology.
Instead, we take advantage of the hypothesis on the Kodaira dimension to run an MMP to reduce to the case when $X$ admits a Mori fiber space structure.
In order to follow this approach, it is necessary to work with singular varieties and pairs.
In fact, we prove a more general version including log canonical 3-fold pairs with nef anti-canonical class.

\begin{theorem} \label{thm: intro_rat_pts_anti-nef}
    Let $p>5$ be a prime number and let $q \coloneqq p^e$ for some $e \geq 1$.
    Let $(X, \Delta)$ be a geometrically integral log canonical projective 3-fold pair over $\mathbb{F}_q$ such that $-(K_X+\Delta)$ is nef. 
	Suppose that the Kodaira dimension of a resolution $\widetilde{X}$ of $X$ is $- \infty$.
    If $q>19$, then $X(\mathbb{F}_q) \neq \emptyset$ holds.
\end{theorem}

We now explain the main ideas in the proof of \autoref{thm: intro_rat_pts_anti-nef}.
A first observation is that the condition that $-(K_X+\Delta)$ is nef is not a birational invariant;
in particular, it is not preserved by the steps of an MMP.
To remedy this, we follow ideas of \cites{FS23,FW21}, and we regard $(X,\Delta)$ as a generalized log Calabi--Yau pair $(X, \Delta, \bM .)$, where $\bM . \coloneqq -(K_X+\Delta)$.
By passing to a partial resolution, we may assume that the Kodaira dimension of $X$ is negative.
Then, by \cite{XZ}, $K_X$ is not pseudo-effective.
By running a suitable MMP, we reach a birational model $Y$
admitting a Mori fiber space structure $Y \to Z$. 
One is then left with two tasks:
to show that (a) the existence of a rational point on $Y$ implies the existence of one on $X$, and (b) to construct a rational point on $Y$.

As for (a), in general, the existence of a rational point is not a birational invariant of singular varieties.
On the other hand, a finite field is a $C_1$-field, and we adapt the proof of \cite{Pie22} to show that the existence of a rational point on $Y$ implies the existence of a rational point on $X$.
To this end, we extend the results of Hogadi--Xu \cite{HX09} on degenerations of Fano varieties in characteristic 0 to the case 3-folds in characteristic $p>5$.
This relies on the systematic use of the extraction of a Koll\'{a}r component and recent vanishing theorems proven in \cite{BK23} (see \autoref{sect: deg_fano}).

For (b), i.e., to construct a rational point on $Y$, we exploit the structure of Mori fiber space $Y \rar Z$.
Concretely, the results in \autoref{sect: deg_fano} allow to lift the existence of a rational point from the base $Z$ to the total space $Y$.
Then, to construct rational points on the base $Z$, one would like to show it is of generalized log Calabi--Yau type.
The statement in characteristic 0 relies on the canonical bundle formula, which is known to fail in characteristic $p>0$. 
Rather than proving that $Z$ is of generalized log Calabi--Yau type, we directly prove the existence of a rational point.
To this end, we utilize various techniques.
In the non-klt case, we construct special models as in \cites{KX16, FS23} to deduce that the base is a variety of log Calabi--Yau type of dimension at most 2; see \autoref{sect: non-klt}.
The existence of rational points on surfaces of generalized log Calabi--Yau type is proven in \autoref{sect: surfaces}, generalizing previous results of Nakkajima \cite{Nak19}.
In the klt case, we use a combination of semi-positivity theorems \cite{Pat14} and the canonical bundle formula in relative dimension 1 \cite{Wit21} (see \autoref{sect: klt-3-folds}).

As a byproduct of this analysis and of the results in \cite{EP23}, we show the surjectivity of the Albanese morphism for 3-folds with nef anti-canonical class.

\begin{theorem}\label{thm: intro_surj_alban}
	Let $k$ be a perfect field of characteristic $p>5$.
	Let $(X,\Delta)$ be a projective log canonical 3-fold pair with $-(K_X+\Delta)$ nef.
	Then, $\alb_X \colon X \to \Alb_X$ is surjective.
\end{theorem}

We conclude the article by studying the existence of rational points on terminal $K$-trivial varieties with non-trivial Albanese morphism.
In the smooth case, this has been addressed in \cite{PZ}*{Corollary 12.9} under the assumption that $X$ is globally $F$-split.
We show that one can drop the assumption that $X$ is globally $F$-split, and we improve the bounds on $q$.
In particular, we obtain the following result.

\begin{theorem} \label{thm: intro_alb_non_trivial}
     Let $p>5$ be a prime number.
     Let $X$ be a geometrically integral canonical $K$-trivial 3-fold over $\mathbb{F}_q$. 
     Suppose that the Albanese morphism $\alb_X$ is not constant.
    If $q>19$, then $X(\mathbb{F}_q) \neq \emptyset$.
\end{theorem}

Here the main difficulty we have to face is that every fiber of $\alb_X$ might be very singular without the globally $F$-split property.
For instance, every fiber may be not log canonical or not even normal.
To overcome this problem, we use Tate's base change formula \cite{PW22} and the connectedness principle \cite{FS23} to show the existence of a rational point on each rational fiber.

\subsection*{Acknowledgements}
We would like to thank A.~Fanelli, G.~Martin, M.~Pieropan, H.~Tanaka, and J.~Waldron for the interesting conversations on the content of this article.
We would like to thank S.~Ejiri and Zs.~Patakfalvi for several useful discussions on their recent work \cite{EP23}.
Lastly, the authors would like to thank the anonymous referee for useful comments that helped the authors improve the content of this work.

\section{Preliminaries}

\subsection{Notation}

\begin{enumerate}
    \item Throughout this article, we fix $p>0$ to be a prime number. 
    The letter $q$ will denote a power $p^e$ for some $e \geq 1$.
    \item Throughout this article, $k$ denotes a field of characteristic $p$.
    We denote by $\overline{k}$ (resp. $k^{\sep}$) an algebraic closure of $k$ (resp. a separable closure).
    \item A $k$-variety (or simply variety) is a quasi-projective integral scheme of finite type over $k$.  
    \item Given a normal $k$-variety $X$ with dualizing complex $\omega_{X}^{\bullet}$, we denote by $\omega_X$ the first non-zero cohomology sheaf of $X$. 
    As $\omega_X$ is a reflexive sheaf of rank 1, there exists a Weil divisor $K_X$ such that $\omega_X \simeq \mathcal{O}_X(K_X).$ The choice of a canonical divisor is fixed throughout the article. 
    Given a proper birational morphism $f \colon Y \to X$ of normal varieties, we choose the canonical divisor $K_Y$ of $Y$ such that $f_*K_Y=K_X$.
    \item Given a normal projective variety $X$ and a $\mathbb{Q}$-Cartier $\mathbb{Q}$-divisor $L$, we denote by $\kappa(L)$ its Iitaka dimension.
    The \emph{Kodaira dimension} $\kappa(X)$ of $X$ is defined to be the Iitaka dimension of the canonical divisor $K_{Y}$ of a smooth projective birational model $Y$.
    \item We say $(X, \Delta)$ is a pair if $X$ is a normal variety, $\Delta$ is an effective $\mathbb{Q}$-divisor on $X$ and $K_X+\Delta$ is $\mathbb{Q}$-Cartier.
    \item For the definition of the singularities of the MMP (as \emph{terminal, klt, dlt, log canonical}), we refer to \cite{kk-singbook}.
    \item Given a projective normal variety $X$, we denote the Albanese morphism by $\alb_X \colon X \to \Alb_X$.  For the properties of the Albanese morphism of projective normal varieties, we refer to \cite{PZ}*{\S~2.6}.
    \item A morphism $f \colon Y \rar X$ of $k$-varieties is a \emph{contraction} if it is proper and $f_*\mathcal{O}_Y=\mathcal{O}_X$.
    \item Given an $\mathbb{F}_p$-scheme $X$, we denote by $F \colon X \to X$ the absolute Frobenius morphism and by $F^{(e)}$ its iterate. 
\end{enumerate}

\subsection{Generalized pairs}

For the language of b-divisors on normal varieties and their properties, we refer to \cite{FW21}*{\S~2.2}.

\begin{definition} \label{def: gen_pair}
    A \emph{generalized sub-pair} $(X, \Delta, \bM .)/Z$ over a variety $Z$ is the datum of 
    \begin{enumerate}
        \item a normal variety $X$, endowed with a projective morphism $X \to Z$; and
        \item a $\mathbb{Q}$-divisor $\Delta$; and
        \item a $\mathbb{Q}$-b-Cartier $\qq$-b-divisor $\bM .$ that is b-nef over $Z$.
    \end{enumerate}
    Moreover, we require that $K_X+\Delta+\bM X.$ is $\mathbb{Q}$-Cartier.
    If $Z=\Spec(k)$, we omit $Z$ from the notation and simply write $(X,\Delta,\bM.)$.
    If $\Delta$ is effective, we say $(X,\Delta,\bM.)$ is a \emph{generalized pair}.
    Furthermore, if $\bM.=0$, we drop $Z$, $\bM.$, and the word generalized from the notation, and we retrieve the usual notions of \emph{sub-pair} and \emph{pair}.
\end{definition}

\begin{definition}
    Let $(X, \Delta, \bM .)/Z$ be a generalized pair and let $\pi \colon X' \to X$ be a projective birational morphism.
    Then, we define $\Delta'$ via the identity
    \[
    K_{X'} +\Delta' +\bM X'. =\pi^*(K_X +\Delta+\bM X.).
    \]
    Given a prime divisor $E$ on $X'$, we define the \emph{generalized log discrepancy} of $E$ with respect to $(X,B,M)/Z$ to be $a_E(X, \Delta, \bM .) \coloneqq 1-\text{coeff}_E(\Delta') $
    .\end{definition}

\begin{definition}
    Let $(X, \Delta, \bM .)/Z$ be a generalized pair. 
    We say that $(X, \Delta, \bM .)/Z$ is \emph{generalized log canonical} (resp. \emph{generalized klt}) if $a_E(X, \Delta, \bM .) \geq 0$ (resp. $a_E(X, \Delta, \bM .) >0$) for all divisors $E$ over $X$.
\end{definition}

\begin{lemma} \label{lem: pertubation-gen}
    Let $(X , \Delta, \bM .)/Z$ be a generalized log canonical pair of dimension at most 3.
    Suppose that $(X,\Delta,\bM.)/Z$ is generalized klt or that $X$ has klt singularities.
    If $A$ is an ample $\mathbb{Q}$-Cartier $\mathbb{Q}$-divisor on $X$, then there exists an effective $\qq$-divisor $\Gamma \sim_{\mathbb{Q},Z} \Delta+\bM X. +A$ such that $(X, \Gamma)$ is klt.
\end{lemma}

\begin{proof}
    First, assume that $(X,\Delta,\bM.)/Z$ is not generalized klt.
    Then, by assumption, $X$ is klt.
    Since $K_X$ and $K_X+\Delta+\bM X.$ are $\qq$-Cartier, then so is $\Delta+\bM X.$.
    Fix $0 < \delta \ll 1$ such that $\delta (\Delta + \bM X.) + A$ is ample over $Z$.
    As $(X,\Delta,\bM.)/Z$ is generalized log canonical and $X$ is klt, we have $(X, (1-\delta) \Delta, (1-\delta) \bM X.)/Z$ is generalized klt and $A' \coloneqq \delta (\Delta + \bM X.) + A$ is ample.
    So, up to replacing $(X,\Delta,\bM.)/Z$ with $(X, (1-\delta) \Delta, (1-\delta) \bM X.)/Z$ and $A$ with $A'$, we can suppose that $(X , \Delta, \bM .)/Z$ is generalized klt.

    Now, let $\pi \colon X' \rar X$ be a projective log resolution where $\bM.$ descends.
    Consider the sub-pair $(X',\Delta')$, where $K_{X'}+\Delta' + \bM X'. = \pi^*(\K X . +\Delta  +\bM X.)$.
    Then, $(X',\Delta')$ is sub-klt.
    Since $\bM X'.$ is nef over $Z$ and $\pi^* A$ is nef and big over $Z$, there exists an effective divisor $F'$ on $X'$ such that, for every $0 < \sigma \ll 1$, we may write $\bM X'. + \pi^* A  \sim_{\qq,Z} H'_\sigma + \sigma F'$, where $H'_\sigma$ is ample (recall that our varieties are quasi-projective).
    We choose $\sigma$ small enough so that $(X',\Delta' + \sigma F')$ is sub-klt; see \cite{FW21}*{Lemma 3.1}.
    Then, we may choose $0 \leq \Gamma' \sim_{\qq} H'_\sigma$ such that $(X',\Delta' + \sigma F' + \Gamma')$ is sub-klt; see \cite{FW21}*{Lemma 3.2}.
    Then, we have that
    \[
    \K X. + \Delta + \bM X. + \epsilon A \sim_{\qq,Z} \K X. + \Delta  + \Sigma,
    \]
    where $\Sigma = \pi_* (\sigma F' + \Gamma')$, and $(X,\Gamma \coloneqq \Delta + \Sigma)$ is a klt pair by construction.
\end{proof}

\begin{lemma} \label{lemma: abundance_big_divisor}
    Let $k$ be a perfect field of characteristic $p > 5$.
    Let $(X, \Delta)$ be a $\mathbb{Q}$-factorial klt pair of dimension at most 3, projective over $Z$. Suppose that $K_X+\Delta$ is nef over $Z$ and $\Delta$ is big over $Z$.
    Then $K_X+\Delta$ is semi-ample over $Z$.
\end{lemma}

\begin{proof}
    As $\Delta$ is big over $Z$, there exist an effective ample divisor $A$ and an effective divisor $E$ such that $\Delta \sim_{\mathbb{Q},Z} A+E$. 
    For $0<\epsilon \ll 1$, the pair $(X, (1-\epsilon)\Delta+\epsilon E)$ is klt .
    Note that $(K_X+\Delta)-(K_X+(1-\epsilon)\Delta+\epsilon E) \sim_{\mathbb{Q},Z} \epsilon A$ is ample over $Z$, and thus we conclude $K_X+\Delta$ is semi-ample over $Z$ by \cite{GNT19}*{Theorem 2.9}.
\end{proof}

The following is an adaptation of \cite{BZ16}*{Lemma 4.4.(2)} that relies on \autoref{lem: pertubation-gen} and \autoref{lemma: abundance_big_divisor}.

\begin{lemma} \label{lemma: bz}
    Let $k$ be a perfect field of characteristic $p > 5$.
    Let $(X,\Delta,\bM.)/Z$ be a $\qq$-factorial generalized klt pair with $\dim(X)=3$.
    Assume that $\K X. + \Delta + \bM X.$ is pseudo-effective over $Z$ and that $\K X. + (1+\alpha)\Delta + (1+\beta) \bM X.$ is big over $Z$ for some $\alpha,\beta \geq 0$.
    Then, we may run a $(\K X. + \Delta + \bM X.)$-MMP over $Z$ with scaling of an ample divisor, and it terminates with a relatively good minimal model.
\end{lemma}

\begin{proof}
    The proof of \cite{BZ16}*{Lemma 4.4.(2)} goes through verbatim by \autoref{lem: pertubation-gen}, the MMP for klt 3-folds \cite{GNT19}*{Theorem 2.13} and \autoref{lemma: abundance_big_divisor}.
\end{proof}

Given a generalized pair $(X,\Delta, \bM.)/Z$, we say a divisor $E$ over $X$ is a \emph{non-klt place} for the generalized pair if $a_E(X,\Delta, \bM.) \leq 0$; in this case, we say that the closed subvariety $\cent_X(E) \subset X$ is a non-klt center for the generalized pair.
Furthermore, if $(X,\Delta,\bM.)/Z$ is generalized log canonical, a non-klt place (resp.~non-klt center) is also called \emph{log canonical place} (resp.~\emph{log canonical center}).
The \emph{non-klt locus} $\nklt(X,\Delta, \bM.)$ is the union of all the non-klt centers of $(X,\Delta, \bM.)/Z$.
We recall the definition of generalized dlt singularities.

\begin{definition} \label{def: gen-dlt-pair}
        A generalized pair $(X,\Delta, \bM .)/Z$ is called \emph{generalized dlt}, if it is generalized log canonical, and, for the generic point $\eta$ of any generalized log canonical center, the following conditions hold:
        \begin{enumerate}
            \item $\bM . = \overline{\bM X.}$ over a neighborhood of $\eta$; and
            \item $(X, \Delta)$ is log smooth in a neighborhood of $\eta$.
        \end{enumerate}
\end{definition}

\begin{definition}
    Let $(X,\Delta, \bM.)/Z$ be a generalized log canonical pair.
    A projective birational morphism $\pi \colon X' \to X$ of normal varieties is a \emph{generalized dlt modification} of $(X,\Delta, \bM.)/Z$ if $X'$ is $\mathbb{Q}$-factorial, every $\pi$-exceptional divisor is a log canonical place of $(X,\Delta,\bM.)/Z$, and, if we write $K_{X'}+\Delta'+\bM X'.  = \pi^*(K_X+\Delta+ \bM X.)$, $(X',\Delta', \bM. )/Z$ is generalized dlt. 
\end{definition}

\begin{proposition} \label{prop:extraction}
    Let $k$ be a perfect field of characteristic $p>5$.
    Let $(X,\Delta, \bM.)/Z$ be a generalized log canonical pair of dimension at most 3.
    Let $\{E_1, \ldots, E_l \}$ be a finite (possibly empty) collection of exceptional divisors over $X$ whose generalized log discrepancy is in $[0,1]$.
    Then, there exists a projective birational morphism $\pi \colon X' \rar X$ from a $\qq$-factorial variety $X'$ such that
    \begin{itemize}
        \item every $E_i$ appears as divisor on $X'$; and
        \item every $\pi$-exceptional divisor is either $E_i$ for some $i$ or is a log canonical place for $(X,\Delta,\bM.)/Z$.
    \end{itemize}
    In particular, we may write $\pi^*(\K X. + \Delta + \bM X.)= \K X'. + \Delta' +\bM X'.$ where $\Delta'$ is an effective divisor.
    Furthermore, $(X',\Delta',\bM .)/Z$ is generalized dlt.
\end{proposition}

\begin{proof}
    Let $f \colon Y \rar X$ be a log resolution where $\bM.$ descends and such that each $E_i$ appears as an exceptional divisor.
    Let $\Sigma_Y= \sum_{i=1}^m E_i$ denote the reduced $f$-exceptional divisor.
    Set $$ \Gamma_Y \coloneqq \Sigma_Y - \sum_{i=1}^l a_{E_i}(X,\Delta,\bM.) E_i.$$ 
    By assumption, we have $a_{E_i}(X,\Delta,\bM.) \in [0,1]$ for all $i$, and each $E_i$ is $f$-exceptional.
    In particular, it follows that $\Gamma_Y$ is effective.
    Since $Y$ is a log resolution and $\bM.$ descends on $Y$, $(Y,\Delta_Y+\Gamma_Y,\bM.)/Z$ is generalized dlt, where $\Delta_Y$ denotes the strict transform of $\Delta$.
    Then, we run a $(\K Y. + \Delta_Y + \Gamma_Y + \bM Y.)$-MMP with scaling over $X$.
    This can be performed, as explained in \cite{FW21}*{proof of Proposition 3.6}.
    
    By construction, we have
    \[
    \K Y. + \Delta_Y + \Gamma_Y + \bM Y.\sim_{\qq,f} \sum_{i=l+1}^m a_{E_i}(X,\Delta,\bM.) E_i = F_Y,
    \]
    where $F_Y$ is an effective, $f$-exceptional divisor.
    More precisely, $\Supp(F_Y)$ coincides with the $f$-exceptional divisors that are neither generalized log canonical places for $(X,\Delta,\bM.)/Z$ nor any of the $E_i$ for $i \in \{1,\ldots,l\}$.
    
    In particular, we have that $\Supp(F_Y) \subset \Supp(\lfloor \Delta_Y + \Gamma_Y \rfloor)$.
    Thus, the MMP terminates by \cite{FW21}*{Proposition 2.12}.
    By the negativity lemma, the only divisors that survive this MMP are $E_1, \dots, E_l$ and the log canonical places extracted.
    Finally, the claim follows as $(Y,\Delta_Y + \Gamma_Y,\bM.)/Z$ is generalized dlt and this property is preserved by a run of an MMP by \cite{Bir19}*{2.13.(2)}.
 \end{proof}

\begin{remark} \label{rmk:extraction}
    Under the assumptions of \autoref{prop:extraction}, further assume that $(X,\Delta,\bM.)/Z$ is generalized klt.
    Then, the set of exceptional divisors for the morphism $\pi$ coincides with $\{E_1,\ldots , E_l\}$.
    Furthermore, if $X$ is $\qq$-factorial, $\pi$ is an isomorphism over the complement of the centers of the $\pi$-exceptional divisors.
\end{remark}

\begin{corollary} \label{prop: gen_dlt_model}
    Let $k$ be a perfect field of characteristic $p>5$.
    Let $(X,\Delta, \bM.)/Z$ be a generalized log canonical pair of dimension at most 3.
    Then, $(X,\Delta, \bM.)/Z$ admits a generalized dlt modification $\pi \colon X' \rar X$.
\end{corollary}

\begin{proof}
    The statement follows immediately from \autoref{prop:extraction} by considering the empty collection of exceptional divisors over $X$.
\end{proof}

\subsection{Varieties of generalized log Calabi--Yau type}

In this subsection, we define varieties of generalized log Calabi--Yau type (cf. \cite{FS23}*{\S~1}) and discuss their basic properties.
We fix a perfect field of characteristic $p>0$.

 \begin{definition} \label{def: gen-logCY}
    A generalized pair $(X, \Delta, \bM.)$ is \emph{generalized log Calabi--Yau} if it is generalized log canonical and $K_X+\Delta + \bM X. \sim_{\mathbb{Q}} 0$.
    If $\bM.$ is 0, we drop the word generalized from the notation.
    We say a projective $k$-variety $X$ is of \emph{generalized log Calabi--Yau type} if there exist an effective $\mathbb{Q}$-divisor $\Delta$ and a b-nef $\qq$-b-Cartier $\qq$-b-divisor $\bM.$ such that $(X,\Delta, \bM.)$ is generalized log Calabi--Yau.
    If $\bM.$ can be chosen to be 0, we say that $X$ is of \emph{log Calabi--Yau type}.
 \end{definition}

\begin{example}
Let $(X,\Delta)$ be a projective log canonical pair with $-(\K X. + \Delta)$ nef.
Then, $X$ can naturally be endowed with the structure of a generalized log Calabi--Yau pair $(X,\Delta,\bM.)$, where we set $\bM. \coloneqq \overline{-(\K X. + \Delta)}$.
\end{example}

\begin{lemma}\label{lem: im-bir-logCY}
Let $(X, \Delta, \bM.)$ be a generalized log Calabi--Yau pair.
Let $\pi \colon X \to Y$ be a birational contraction of normal projective varieties.
Then $(Y, \pi_*\Delta, \bM.)$ is a generalized log Calabi--Yau pair. 

In particular, given a projective birational contraction $\pi \colon X \to Y$, if $X$ is of (generalized) log Calabi--Yau type, then so is $Y$. 
\end{lemma}

\begin{proof}
    As we have $K_X+\Delta + \bM X. \sim_{\mathbb{Q}} 0$, it follows that $K_Y+\pi_*\Delta + \bM Y.$ is $\mathbb{Q}$-Cartier and $\mathbb{Q}$-linearly equivalent to zero.
    By the negativity lemma \cite{7authors}*{Lemma 2.16}, we have  
    $K_X+\Delta + \bM X. = \pi^*(K_Y+\pi_*\Delta  +\bM Y.)$, concluding that $(Y, \pi_*\Delta, \bM .)$ has generalized log canonical singularities.
\end{proof}

\begin{lemma} \label{lem: descent-finite-morph}
    Let $k$ be a field of characteristic $p > 0$.
    Let $f \colon X \to Y$ be a generically finite projective morphism of degree $\deg(f)<p$ between normal varieties over $k$.
    If $X$ is of (generalized) log Calabi--Yau type, then so is $Y$.
\end{lemma}

\begin{proof}
    If $X \to X' \to Y$ is the Stein factorization of $f$, then $X'$ is of generalized log Calabi--Yau type by \autoref{lem: im-bir-logCY} and we can reduce to the case where $f$ is finite.
    Note that $f$ is separable as $p$ does not divide $\deg(f)$.
    The proof of \cite{HL21}*{Theorem 4.5} applies verbatim in our case as 
    the hypothesis on $\deg(f)<p$ 
    permits to apply the Hurwitz formula \cite{kk-singbook}*{2.41} as $f$ is tame along every prime divisor.
\end{proof}

\begin{lemma} \label{lemma:perturbation-global}
    Let $(X,\Delta,\bM.)$ be a generalized log Calabi--Yau pair.
    Assume that $X$ is a klt variety.
    Further assume that $X$ admits a projective contraction $f \colon X \rar Z$ with $-K_X$ ample over $Z$, and fix an ample $\qq$-Cartier $\qq$-divisor $H$ on $Z$.
    Then, for every $\epsilon > 0$, we may find a boundary $\Gamma$ such that $(X,\Gamma)$ is klt and $K_X + \Gamma \sim_{\qq} \epsilon f^*H$. 
\end{lemma}

\begin{proof}
    By assumption, $\Delta+ \bM X. \sim_{\qq} -K_X$ is $\qq$-Cartier and $f$-ample.
    In particular, $\Delta+ \bM X.$ is $f$-ample.
    Fix $0 < \delta \ll \epsilon$ such that $\delta(\Delta + \bM X.)+\epsilon f^*H$ is ample.
    Then, we conclude by applying \autoref{lem: pertubation-gen} to the generalized pair $(X,(1-\delta)\Delta,(1-\delta)\bM.)$ and the ample divisor $\delta(\Delta + \bM X.)+\epsilon f^*H$.
\end{proof}

\section{Rational points on surfaces of generalized log Calabi--Yau type} \label{sect: surfaces}

In this section, $p$ is a prime number and $q=p^e$ for some integer $e \geq 1$.
We show the existence of rational points for surface pairs with nef anti-canonical class defined over a finite field $\mathbb{F}_q$, with $q>19$. 
These results should be compared with the ones obtained by Nakkajima in \cite{Nak19}.
We start with the case of curves.

\begin{lemma}\label{lem: rat_pts_curves}
    Let $C$ be a geometrically integral regular projective curve with $-K_C$ nef over $\mathbb{F}_q$ (i.e., $g(C) \leq 1$).
    Then, we have $C(\mathbb{F}_q) \neq \emptyset$. 
\end{lemma}

\begin{proof}
    If $g(C)=0$, $C$ is a conic by \cite{kk-singbook}*{Lemma 10.6} and we apply the Chevalley--Warning's theorem \cite{Ser73}*{\S~I.I.2 Theorem 3}. 
    If $g(C)=1$, the base change $C_{\overline{\mathbb{F}_q}}$ is an elliptic curve and thus we conclude by \cite{Lan55}*{Theorem 3}.
\end{proof}

Now, we show that the image of a surface of generalized log Calabi--Yau type via a contraction is a variety of generalized log Calabi--Yau type.

\begin{proposition} \label{prop:im_gen_CY_Surfaces}
    Let $k$ be a perfect field of characteristic $p>0$ and let $(X, \Delta, \bM .)$ be a generalized log Calabi--Yau surface pair over $k$.
    If $f \colon X \to Y$ is a contraction, then $Y$ is of generalized log Calabi--Yau type.
\end{proposition}

\begin{proof}
    If $\dim(Y)=0$, the claim is immediate.
    Similarly, if $\dim(Y)=2$, $f$ is birational and the claim follows from \autoref{lem: im-bir-logCY}.
    Thus, it remains to prove the case where $\dim(Y)=1$.
    Equivalently, we have to show that the geometric genus $g(Y)$ is at most $1$. 
    For this statement, as $k$ is perfect, we may assume that it is algebraically closed.

    By passing to the minimal resolution, we can suppose $X$ is smooth.
    Then, we run a $K_X$-MMP over $Y$.
    By replacing $X$ with the outcome of this MMP and $f$ with the induced morphism to $Y$, we may assume that either $K_X$ is $f$-nef or $f$ is a $K_X$-Mori fiber space.

    Suppose that $K_X$ is $f$-nef.
    Since $-K_X$ is pseudo-effective, we conclude $K_X \equiv_f 0$.
    By abundance for surfaces \cite{Tan14}*{Theorem 6.7}, $K_X$ is $f$-semiample and thus there exists a $\mathbb{Q}$-divisor $L$ on $Y$ such that $K_X \sim_{\mathbb{Q}} f^*L$. 
    As $-f^*L \sim_{\mathbb{Q}} \Delta+\bM .$ is pseudo-effective and $Y$ is a curve, we deduce that $-L$ is nef.
    Therefore $-K_X$ is nef. 
    In particular $\kappa(X)$ is either 0 or $-\infty$ and thus we conclude that the genus of $Y$ is at most 1 by the Iitaka's conjecture \cite{CZ15}*{Theorem 1.3}.

    Suppose $f \colon X \to Y$ is a $K_X$-Mori fiber space. Then, $\Delta+\bM X.$ is $f$-ample. 
    Fix a positive integer $n$ and an ample $\qq$-Cartier $\qq$-divisor $H$ on $C$.
    By \autoref{lemma:perturbation-global}, we may find a boundary $\Gamma_n$ such that $(X,\Gamma_n)$ is klt and $K_X + \Gamma_n \sim_\qq \frac 1 n f^* H$.

    Suppose by contradiction that $g(Y) \geq 2$. 
    As $k$ is a perfect field, the generic fiber $X_\eta$ is geometrically reduced by \cite{Sch10}. As $X_\eta$ is a geometrically reduced and it is a conic by \cite{kk-singbook}*{Lemma 10.6}, it is smooth. 
    Therefore, by \cite{Wit21}*{Theorem 3.4}, there exists a purely inseparable morphism $\pi_n \colon T_n \to Y$ such that ${\pi}^*_n \frac 1 n H \sim_\qq t_n K_{T_n}+(1-t_n) \pi^* K_Y+E_n$ for some $t_n \in [0,1]$ and $E_n \geq 0$.
    By the proof of \autoref{lemma:perturbation-global}, $\Gamma_n$ is a small perturbation of $\Delta$.

    \begin{claim}\label{claim:surfaces}
        Up to considering a subsequence, the morphism $\pi_n \colon T_n \rar Y$ does not depend on $n$, $\{E_n\}_{n \in \mathbb N}$ is convergent up to $\qq$-linear equivalence, and $\{ t_n \}_{n \in \mathbb N}$ is a converging sequence with limit $t_{\infty} \in [0,1]$.
    \end{claim}

    The setup of \autoref{claim:surfaces} is completely analogous to the one of \autoref{claim:epsilon}.
    We observe that, unlike in the case of \autoref{claim:epsilon}, here a small perturbation of $\Delta$ is in place due to the use of \autoref{lemma:perturbation-global}.
    Yet, as the perturbation is small, the same proof as in \autoref{claim:epsilon} goes through, with the due bookkeeping of the small perturbation coming from \autoref{lemma:perturbation-global}.
    Thus, for the details, we refer to the proof of \autoref{claim:epsilon}.
    
    Hence, by \autoref{claim:surfaces}, we may assume that $T=T_n$ and $\pi=\pi_n$ are independent of $n$, $\{t_n\}_{n \in \mathbb N}$ is a converging sequence with limit $t_\infty \in [0,1]$, and $\{E_n\}_{n \in \mathbb N}$ is convergent up to $\mathbb Q$-linear equivalence.
    Then, considering the limit as $n \to \infty$ and arguing as in the proof of \autoref{prop: point_base_conic_bundle_gen}, we conclude that $-K_T$ is pseudo-effective.
    Thus, we conclude that $g(T) \leq 1$.
    Since $\pi$ is a universal homeomorphism, we obtain that $g(Y) \leq 1$, leading to the sought contradiction.
\end{proof}

For the case of $K$-trivial surfaces with canonical singularities, we use the boundedness of their $l$-adic Betti number and the Weil conjectures.
We thank G. Martin for suggesting the proof in the case of Enriques surfaces.

\begin{proposition} \label{prop: rat_pts_reg_K_triv}
    Let $X$ be a geometrically integral projective surface over $\mathbb{F}_q$ with canonical singularities and numerically trivial canonical class.
    If $q > 19$, then we have $X(\mathbb{F}_q) \neq \emptyset$.
    More precisely, if $X'$ denotes the minimal resolution of $X$, we have that:
    \begin{itemize}
        \item if $X=X'$ is geometrically an Abelian surface or a bielliptic surface, we may take $q \geq 2$;
        \item if $X'$ is geometrically a K3 surface, we may take $q >19$;
        \item if $X'$ is geometrically an Enriques surface, we may take $q \geq 2$.
    \end{itemize}
\end{proposition}

\begin{proof}
    Replacing $X$ with its minimal resolution, we can suppose $X$ is smooth.
    As $\mathbb{F}_q$ is perfect, we can apply the classification of Bombieri--Mumford \cite{BM77} on $X_{\overline{{\mathbb F}_q}}$.
    We proceed by cases according to the classification.
    
    If $X_{\overline{{\mathbb F}_q}}$ is an Abelian surface, then $X$ has a rational point by \cite{Lan55}.
    
    If $X_{\overline{{\mathbb F}_q}}$ is bielliptic, we consider the Albanese morphism $X \to \Alb_X$.
    This is a surjective fibration that is isotrivial in the \'etale topology, and $\Alb_X$ is 1-dimensional by \cite{BM77}.
    Therefore $\Alb_X(\mathbb{F}_q) \neq \emptyset$ by \autoref{lem: rat_pts_curves}.
    Fix $c \in \Alb_X$ a rational point. 
    As the fiber over a closed point is either a geometrically irreducible smooth curve of genus 1 or a geometrically irreducible cuspidal curve of genus 1, we conclude that $X_c$ has a rational point by \autoref{lem: rat_pts_curves}.

    Now, let $l$ be a prime number different from $p$.
    If $X_{\overline{{\mathbb F}_q}}$ is a K3 surface, the $l$-adic Betti numbers are $b_1(X)=b_3(X)=0$ and $b_2(X) = 22$ 
    by \cite{BM77}.
    Using the Lefschetz trace formula and the estimates on the eigenvalues of the Frobenius action $F^e$ on the \'etale cohomology  provided by the Weil conjectures \cite{Del74}*{Théorème 1.6}, we deduce for $l \neq p$ that
    $$|X(\mathbb{F}_q)-1-q^2|= |\Tr((F^e)|{H^2_{et}(X, \mathbb{Q}_l})| \leq b_2(X) \cdot q, $$ 
    thus concluding.

    If $X$ is an Enriques surface, the crystalline cohomology groups of $X$ are $H^0_{\text{crys}}(X/K)=H^4_{\text{crys}}(X/K)=K$, $H^2_{\text{crys}}(X/K)=K^{\oplus 10}$ and $H^1_{\text{crys}}(X/K)=H^3_{\text{crys}}(X/K)=0$ by \cite{Ill79}*{Proposition 7.3.5}, where $K$ is the fraction field of the ring of Witt vectors $W(\mathbb{F}_q)$.
    By \cite{Ill79}*{Proposition 7.3.6}, the slopes of $H^2_{\text{crys}}(X/K)$ are all 1.
    Therefore, by the Lefschetz trace formula in crystalline cohomology \cite{Ber74}*{Chapitre VII.3}, we conclude
    $$
    |X(\mathbb{F}_q)| = \sum (-1)^i \Tr( F^e|{H^{i}_{\text{crys}}(X/K)}) = 1 \text{ mod } q .
    $$
\end{proof}

Now, we treat the case of surfaces of generalized log Calabi--Yau type.

\begin{proposition} \label{prop:rat_pts_CY?surf}
Let $X$ be a geometrically integral surface of generalized log Calabi--Yau type over $\mathbb{F}_q$. 
Then either
\begin{enumerate}
    \item the canonical divisor of the minimal resolution is not pseudo-effective, and $X(\mathbb{F}_q) \neq \emptyset$; or
    \item $X$ has canonical singularities, $K_X \equiv 0$ and $X(\mathbb{F}_{q}) \neq \emptyset$ for $q>19$.
\end{enumerate}
\end{proposition}

\begin{proof}
    By \autoref{prop:extraction}, we can replace $X$ with its minimal resolution.
    In case (b), we conclude by \autoref{prop: rat_pts_reg_K_triv}.
    In case (a), either $\Delta>0$ or $\bM X.$ is not numerically trivial with $K_X+\Delta + \bM X. \sim_{\mathbb{Q}} 0$.
    As $K_X$ is not pseudo-effective, a $K_X$-MMP will end with a smooth surface $Y$ endowed with a Mori fiber space structure $Y \to Z$.
    By \autoref{lem: im-bir-logCY}, $Y$ is of generalized log Calabi--Yau type.
    If $\dim(Z)=0$, then $Y$ is 
    a del Pezzo surface and $Y(\mathbb{F}_q) \neq \emptyset$ by  \cite{Esn03}*{Corollary 1.3}.
    If $\dim(Z)=1$, we have that $g(Z) \leq 1$ by \autoref{prop:im_gen_CY_Surfaces} and thus $Z(\mathbb{F}_q) \neq \emptyset $ by \autoref{lem: rat_pts_curves}.
    Let $z \in Z$ be a $\mathbb{F}_q$-point. 
    As the fiber $F_z$ is a conic with $H^0(F_z, \mathcal{O}_z)=\mathbb{F}_q$ by \cite{BT22}*{Proposition 2.18}, we conclude it has a rational point by the Chevalley--Warning's theorem \cite{Ser73}*{\S~I.I.2 Theorem 3}.
    Since $Y$ is smooth, the existence of an $\mathbb F _q$-point on $X$ follows from the Lang--Nishimura's lemma \cite{BT22}*{Lemma 6.11}.
\end{proof}

\begin{remark}
Consider the Fermat quartic surface $X \coloneqq \left\{x_0^4+x_1^4+x_2^4+x_3^4=0 \right\} \subset \mathbb{P}^3_{\mathbb{F}_5}$. 
By Fermat's little theorem, we have that $a^4 \in \left\{0,1 \right\}$ for  every $a \in \mathbb{F}_5$ and this shows that $X(\mathbb{F}_5) =\emptyset$.
We do not know if $q>19$ is sharp for the existence of rational points on K3 surfaces over finite fields. 
\end{remark}

\begin{corollary} \label{cor: surf_bir_rat_pts}
    Let $X$ be a geometrically integral projective surface over $\mathbb{F}_q$ that is birational to a surface $Y$ of generalized log Calabi--Yau type.
    If $q>19$, then $X(\mathbb{F}_{q}) \neq \emptyset$.
\end{corollary}

\begin{proof}
    By passing to the minimal resolution, we can suppose $Y$ is smooth.
    Then $Y(\mathbb{F}_{q}) \neq \emptyset$ by \autoref{prop:rat_pts_CY?surf} and we conclude by the Lang--Nishimura's lemma \cite{BT22}*{Lemma 6.11}.
\end{proof}

\subsection{Surjectivity of the Albanese morphism}

We now show the surjectivity of the Albanese morphism for surfaces of generalized log Calabi--Yau type.

\begin{lemma} \label{lem: im_alb}
    Let $k$ be a perfect field, and let $f \colon Y \to X$ be a surjective morphism of normal projective varieties. 
    If $\alb_Y$ is surjective, then so is $\alb_X$.
\end{lemma}

\begin{proof}
    Consider the commutative diagram obtained by the universal properties of the Albanese morphism:
      $$
    \xymatrix{
    Y \ar@{->>}[r]^{\alb_Y} \ar[d] & \Alb_Y \ar[d]^{}   \\
    X \ar[r]^{\alb_X} & \Alb_X.
    }
    $$
    As the image of $\alb_X$ generates $\Alb_X$ as a group variety and the map $\Alb_Y \to \Alb_X$ is a group homomorphism, we conclude that $\Alb_Y \to \Alb_X$ is surjective and thus we deduce that also $\alb_X$ is surjective.
\end{proof}

\begin{proposition} \label{prop: alb_surf}
Let $k$ be a perfect field of characteristic $p>0$ and let $X$ be a geometrically integral surface of generalized log Calabi--Yau type over $k$. Then $\alb_X$ is surjective.
\end{proposition}

\begin{proof}
    We can suppose $k$ is algebraically closed.
    By \autoref{lem: im_alb}, we may pass to the minimal resolution and assume that $X$ is smooth.
    By running a $K_X$-MMP, we either end up with a Mori fiber space structure $f \colon Y \to B$ or with a smooth surface $Y$ with trivial canonical class.
    In the first case, as the fibers of $f$ are rationally connected, we deduce that $\alb_X$ is surjective if and only $\alb_B$ is surjective.
    If $B$ is a point, it is immediate, while if $B$ is a curve, we conclude by \autoref{prop:im_gen_CY_Surfaces}.
    In the second case, we can conclude by the classification of Bombieri--Mumford \cite{BM77} (or the recent result \cite{EP23}*{Corollary 1.4}).
\end{proof}
\section{Degenerations of varieties of Fano type} \label{sect: deg_fano}

In this section, we extend the results of \cite{HX09} to 3-folds admitting a contraction of Fano type over a perfect field of characteristic $p>5$. 
In view of the birational MMP for 3-folds developed by \cite{HW22}, the restriction on the characteristic is due to the use of the Kawamata--Viehweg vanishing theorem for del Pezzo surfaces of \cite{ABL22} and their applications \cites{HW19, BK23}. 
Our proof relies on a systematic use of the extraction of a Koll\'ar component as developed in \cite{GNT19}.

\begin{definition} \label{def: contraction_Fano_type}
    Let $f \colon X \to Z$ be a projective contraction of normal varieties.
    We say $f$ is a \emph{contraction of Fano type} if there exists $\Delta \geq 0$ such that $(X, \Delta)$ has klt singularities and $-(K_X+\Delta)$ is $f$-big and $f$-nef.
    If $Z=\Spec(k)$, we say $X$ is a \emph{variety of Fano type}.
    If $\dim(X)=2$ and $Z=\Spec(k)$, we say that $X$ is of \emph{del Pezzo type}.
\end{definition}

We start with the analog of \cite{HX09}*{Theorem 1.3}.

\begin{proposition} \label{prop: Hogadi-Xu-deg}
Let $k$ be a perfect field of characteristic $p>5$.
Let $(X, \Delta)$ be a geometrically integral klt 3-fold pair over $k$ and let $x$ be a closed point of $X$.
Then, there exists a log resolution $\pi \colon Y \to X$ such that 
\begin{enumerate}
    \item there is an irreducible $\pi$-exceptional divisor such that $\cent_X(E)=x$;
    \item $E$ is a geometrically irreducible $k(x)$-variety; and
    \item $E$ is a separably rationally connected surface.
\end{enumerate}
\end{proposition}

\begin{proof}
    Since $(X,\Delta)$ is klt, it is of Fano type over $Z=X$.
    Let $g \colon X' \to X$ be a projective birational morphism extracting a Koll\'{a}r component $E$ over $x$ as in \cite{GNT19}*{Proposition 2.15}.
    Notice that, up to shrinking $X$ around $x$, by \cite{GNT19}*{Proposition 2.15.(4)}, we may assume that $g$ is an isomorphism over $X \setminus \{x\}$.
    Since $(X,\Delta)$ is klt, we have $-E \sim_{\mathbb Q , g} K_{X'} + g_*^{-1} \Delta + (a-1) E$ for some $a<1$.
    As $E$ is a normal surface of del Pezzo type and $-E$ is $\pi$-big and $\pi$-nef, we can apply \cite{BK23}*{Proposition 2.3} together with \cite{ABL22}*{Theorem 1.1} to show that $R^1g_* \mathcal{O}_{X'}(-E)=0$.
    Thus, $g_* \mathcal{O}_{X'} \simeq \mathcal{O}_X \to g_* \mathcal{O}_E$ is surjective.
    This implies that $E$ is geometrically irreducible over $k(x)$.
    
     As $E$ is a normal surface of del Pezzo type over a perfect field, we have $E_{\overline{k}}$ is a rational surface by \cite{Tan-x-meth}*{Fact 3.4 and Theorem 3.5} and thus $E$ is separably rationally connected by \cite{Kol96}*{IV.3.2.5}.
     
    Let $\psi \colon Y \to X'$ be a log resolution for $(X', E)$. As $E$ is normal, its strict transform $E_{Y} \coloneqq \psi_*^{-1}(E)$ is also geometrically irreducible and separably rationally connected and thus we conclude by setting $\pi \coloneqq g \circ \psi$.
\end{proof}

The following is a generalization of the Lang--Nishimura's lemma for klt 3-folds over finite fields. 
A similar result for fields of characteristic 0 is contained in \cite{Pie22}*{Theorem 1.3}.

\begin{corollary} \label{cor:rat-pts-finite-fields}
Let $p$ be a prime number with $p >5$, and let $q \coloneqq p^e$ for some integer $e \geq 1$.
Let $(X, \Delta)$ be a geometrically integral klt 3-fold pair over the finite field $\mathbb{F}_q$ and let $\psi \colon X \dashrightarrow X'$ be a birational morphism of projective varieties.
If $X(\mathbb{F}_q) \neq \emptyset$ holds, then $X'(\mathbb{F}_q) \neq \emptyset$ holds.
\end{corollary}

\begin{proof}
    Let $x \in X$ be an $\mathbb{F}_q$-rational point.
    Let $\pi \colon Y \to X$ be a log resolution that factors through the birational model given by \autoref{prop: Hogadi-Xu-deg}, and let $E$ denote the geometrically irreducible separably rationally connected exceptional divisor constructed in \autoref{prop: Hogadi-Xu-deg}. 
    Without loss of generality, we can suppose $Y$ admits a birational contraction onto $X'$. 
    By \cite{Esn03}*{Corollary 1.3}, we have $E(\mathbb{F}_q) \neq \emptyset$ and this implies $X'(\mathbb{F}_q) \neq \emptyset$. 
\end{proof}

The following is the analog of \cite{HX09}*{Theorem 1.2} for 3-folds admitting a fibration of Fano type in characteristic $p>5$.

\begin{theorem} \label{thm: HX-deg-charp}
    Let $k$ be a perfect field of characteristic $p>5$ and
    let $f \colon X \to Z$ be a projective contraction of quasi-projective varieties, where $\dim(X)=3$ and $\dim(Z) \geq 1$.
    Suppose $f$ is a contraction of Fano type.
    
    If $z \in Z$ is a closed point, then the fiber $X_z$ contains a separably rationally connected subvariety that is geometrically irreducible over $k(z)$.
\end{theorem}

\begin{proof}
    By \cite{GNT19}*{Proposition 2.15}, we may find a commutative diagram of normal quasi-projective varieties
    $$
    \xymatrix{
    W \ar[r]^{\psi} \ar^{\varphi}[d] & Y \ar[d]^{g}   \\
    X \ar[r]^f & Z
    }
    $$
    where  $W$ is smooth, $\psi$ and $\varphi$ are projective birational morphism,
    $(Y, \Delta_Y)$ is plt, $-(K_Y+\Delta_Y)$ is ample over $Z$, and $E \coloneqq (g^{-1}(z))_{\text{red}}=\lfloor{\Delta_Y \rfloor}$ is a normal surface of del Pezzo type such that $-E$ is nef over $Z$.
    By \cite{BK23}*{Proposition 2.3} together with \cite{ABL22}*{Theorem 1.1}, we deduce that $R^{i}g_* \mathcal{O}_Y(-E)=0$ for $i>0$
    in a neighborhood of $z$.
    In particular, $g_*\mathcal{O}_Y \simeq \mathcal{O}_Z \to g_*\mathcal{O}_E $ is surjective and thus $E$ is geometrically irreducible.
    As $E$ is a surface of del Pezzo type, it is geometrically rational by \cite{NT20}*{Proposition 2.26}.

    As $E$ is normal, we have that the strict transform $E_W$ is geometrically irreducible over $k(z)$ and therefore also $\text{center}_X(E_W)$. 
    Finally, as $E_W$ is separably rationally connected by \cite{Kol96}*{IV.3.2.5} we conclude that so is $\cent_X(E_W)$ and this concludes the proof.
\end{proof}

The following gives an example where an irreducible fiber of a Mori fiber space over a perfect field is not necessarily geometrically irreducible.
In particular, in the statement of \autoref{thm: HX-deg-charp}, the separably rationally connected subvariety of the fiber $X_z$ need not be an irreducible component of $X_z$.

\begin{example}
Let $k$ be a field such that there exists a closed point $x \in \mathbb{P}^2_{k}$ with the following property:
The residue field $k(x)$ is a Galois extension of degree 4 of $k$ and the 4 closed points $x_{\overline{k}}=\left\{ x_1, x_2, x_3, x_4 \right\} \subset \mathbb{P}^2_{\overline{k}}$ are in general position.
By the general position hypothesis, the blow-up $X \to \mathbb{P}^2_k$ at $x$ is a del Pezzo surface of degree 5.
The pencil $V$ of conics passing through $x$ becomes a basepoint free pencil on $X$ inducing a Mori fiber space $f \colon X \to \mathbb{P}^1_k$ of relative dimension 1. 
The discriminant locus of $f$ is the disjoint union of 3 rational points on $\mathbb{P}^1_k$ and the singular fibers are integral $k$-conics with a unique rational point (the singular point of the fiber).
Note that the singular fibers are not geometrically integral and not rational.
\end{example}

We recall a consequence of \cite{Esn03} for rationally connected varieties we will use in combination with \autoref{prop: Hogadi-Xu-deg} and \autoref{thm: HX-deg-charp} to construct rational points over finite fields.

\begin{proposition} \label{prop: strenghten-Esnault}
	Let $X$ be a geometrically integral rationally connected projective variety over a finite field $\mathbb{F}_q$.
	If we assume $\dim(X) \leq 3$, then $X(\mathbb{F}_q) \neq \emptyset$ holds.
\end{proposition}

\begin{proof}
	Let $X' \to X$ be a resolution of singularities, whose existence is guaranteed by the hypothesis on the dimension \cite{CP19}.
	As $X'$ is also rationally connected by \cite{Kol96}*{Proposition IV.3.3.3}, we have $X'(\mathbb{F}_q) \neq \emptyset$ by \cite{Esn03}*{Corollary 1.3}, thus concluding. 
\end{proof}
\section{Rational points on 3-folds of generalized log Calabi--Yau type}

In this section, we study the existence of rational points on 3-folds admitting a generalized log Calabi--Yau structure over finite fields.

In the following, the restriction on the characteristic of the base field is required for using the MMP for 3-folds over perfect fields \cite{GNT19}, the results of \autoref{sect: deg_fano} based on the vanishing theorems of \cite{BK23}, and to control the geometric generic fiber of 3-fold Mori fiber spaces \cites{PW22, BT22}.

\subsection{Non-klt generalized log Calabi--Yau 3-fold pairs} \label{sect: non-klt}

We prove the following result on birational models of Calabi--Yau 3-folds in positive characteristic, extending \cite{KX16}*{Corollary 59}.
We will mostly follow \cite{FS23}*{Theorem 4.2}.
As we work in positive characteristic, the canonical bundle formula fails, in general, see \cite{Tan20}, so we need to do some work to prove that the base $Z$ of a special model thus constructed is of generalized log Calabi--Yau type. 

\begin{definition}
    Let $F$ be an effective $\mathbb R$-divisor on a variety $X$, and let $D$ be an $\mathbb R$-divisor on $X$.
    We say that $F$ \emph{fully supports} $D$ if $D$ is effective and the support of $D$ coincides with the support of $F$.
    Moreover, given a proper morphism of varieties $f \colon X \rightarrow Z$, we say that $F$ fully supports an $f$-ample divisor if $F$ fully supports an effective $f$-ample $\mathbb R$-divisor $H$.
\end{definition}

The following is a technical generalization of \cite{FW21}*{Proposition 2.15}.

\begin{lemma} \label{lemma:MMP_not_pseff}
Let $k$ be a perfect field of characteristic $p > 5$.
Let $(X,\Delta,\bM.)/Z$ be a
$\qq$-factorial
generalized dlt pair with $\dim (X) \leq 3$ defined over $k$.
Let $f \colon X \rar Y$ be a contraction over $Z$.
Assume that $(\K X. + \Delta + \bM X.)$ is not pseudo-effective over $Y$.
Then, we may run a $(\K X. + \Delta + \bM X.)$-MMP over $Y$, which terminates with a Mori fiber space.
\end{lemma}

\begin{proof}
We follow the proof of \cite{FW21}*{Proposition 2.15}.
Recall that our varieties are quasi-projective, and let $A$ be an ample divisor on $X$ such that $K_X + \Delta + \bM X. + A$ is $f$-nef.
We claim that we can run a $(K_X + \Delta + \bM X. )$-MMP over $Y$ with scaling of $A$.
Since $(K_X + \Delta + \bM X. )$ is not pseudo-effective over $Y$, this MMP is the same as a $(K_X + \Delta + \bM X. + \epsilon A)$-MMP with scaling of $(1-\epsilon)A$, where $0 < \epsilon \ll 1$.
By \autoref{lem: pertubation-gen}, there exists an effective $\qq$-divisor $\Gamma$ such that $\Gamma \sim_{\qq,Y} \Delta + \bM X. + \epsilon A$ and $(X,\Gamma)$ is klt.
Then, the sought MMP exists and terminates with a Mori fiber space by \cite{GNT19}*{Theorem 2.13}. 
\end{proof}

In \cite{FW21}, the connectedness principle is only stated for pairs with anti-nef log canonical divisor.
Nevertheless, the arguments in \emph{op.~cit.} go through also in the case of generalized pairs of log Calabi--Yau type, by applying the slightly more general perturbations arguments in \autoref{lem: pertubation-gen} and their immediate consequences as \autoref{lemma:MMP_not_pseff}.
In particular, we obtain the following statement.
For the precise terminology involved, we refer to \cite{FW21}.

\begin{corollary} \label{conn_principle}
    Let $k$ be a perfect field of characteristic $p>5$.
    Let $f \colon X \rar Z$ be a projective morphism between normal varieties that are quasi-projective over $\Spec(k)$.
    Assume that $\dim (X) \leq 3$ and that $(X,\Delta,\bM.)/Z$ is a generalized pair with $\K X. + \Delta + \bM X. \sim_{\qq,Z} 0$.
    Fix $z \in Z$, and assume $f^{-1}(z)$ is connected as $k(z)$-scheme.
    Then, if $f \sups -1. (z) \cap \nklt(X,\Delta,\bM.)$ is disconnected as $k(z)$-scheme, the following holds:
	\begin{itemize}
		\item[(1)] $(X,\Delta,\bM.)$ is generalized log canonical in a neighborhood of $f^{-1} (z)$;
	    \item[(2)] $f^{-1}(z)\cap \nklt(X,\Delta,\bM.)$ has two connected components, and, {after an \'etale base change, $\nklt(X,\Delta,\bM.)$ has two connected components over $z$, and each of the connected components of $f^{-1}(z)\cap \nklt(X,\Delta,\bM.)$ corresponds to one of the connected components of $\nklt(X,\Delta,\bM.)$}; and
		\item[(3)] there is an \'etale morphism $(z' \in Z') \rar (z \in Z)$ and a projective morphism $T' \rar Z'$ such that $k(z')=k(z)$ and the crepant pull-back of $(X,\Delta,\bM.)$ to $X \times_Z Z'$ is birational to a standard $\pr 1.$-link over $T'$.
	\end{itemize}
\end{corollary}

\begin{proof}
    The proof goes through verbatim as in \cite{FW21}*{Theorem 1.2}.
    Whenever the initial pair $(X,\Delta)$ is replaced by a dlt modification, we replace the generalized pair $(X,\Delta,\bM.)/Z$ with a generalized dlt modification.
    Then, \cite{FW21}*{Proposition 2.15} is superseded by the slightly more general \autoref{lemma:MMP_not_pseff}.
    Similarly, the MMP needed in the proof of \cite{FW21}*{Proposition 3.6} can be run for a generalized dlt pair by the perturbation tools guaranteed by \autoref{lem: pertubation-gen}.
\end{proof}

We collect a few lemmas we use multiple times in the proof of the existence of special birational model of log Calabi--Yau 3-fold pairs.

\begin{lemma}\label{lemma:rel_ample}
    Let $(X,\Delta,\bM.)/Z$ be a generalized log canonical pair that is $\qq$-factorial.
    Further assume that $\mathrm{Nklt}(X,\Delta,\bM.) = \Supp(\Delta^{=1})$.
    Let $\pi \colon X' \to X$ be a birational contraction of normal varieties that only extracts generalized log canonical places of $(X,\Delta,\bM.)/Z$.
    Write $ K_{X'}+\Delta' +\bM_{X'}. \sim_{\mathbb{Q}} \pi^*(K_X+\Delta+\bM_{X}.)$.
    If $\Delta^{=1}$ fully supports a $\qq$-divisor that is ample over $Z$, then so does $(\Delta')^{=1}$.
\end{lemma}

\begin{proof}
    Since $X$ is $\qq$-factorial, there exists a $\pi$-ample divisor $E' \leq 0$ that is fully supported on the exceptional locus of $\pi$.
    Let $H$ be a $\qq$-divisor with $\Supp(H)=\Supp(\Delta^{=1})$ that is ample over $Z$.
    By assumption, the centers in $X$ of $\pi$-exceptional divisors are contained in $\Supp(H)$.
    Therefore, $\Supp(\pi^*H)$ contains all the $\pi$-exceptional divisors.
    In particular, $\pi^*H + \epsilon E'$ is effective and ample over $Z$ for $0 < \epsilon \ll 1$.
\end{proof}

\begin{lemma}\label{lemma:extraction}
    Let $k$ be a perfect field of characteristic $p > 5$.
    Let $(X,\Delta,\bM.)/Z$ be a $\mathbb{Q}$-factorial generalized log canonical pair over $\Spec(k)$ with $\dim(X) \leq 3$.
    Assume that $X$ is quasi-projective and $\qq$-factorial and also admits a structure of generalized klt pair $(X,\Delta',\mathbf N)/Z$.
    Let $E$ be a generalized log canonical place for $(X,\Delta,\bM.)/Z$ exceptional over $X$.
    Then, there exists an extremal extraction $\pi \colon X' \to X$ that extracts exactly $E$.
\end{lemma}

\begin{proof}
    By assumption and the linearity of log discrepancies in convex combinations, for $0 < \epsilon \ll 1$, $(X,(1-\epsilon)\Delta + \epsilon \Delta',(1-\epsilon)\bM. + \epsilon \mathbf N)/Z$ is generalized klt and the generalized log discrepancy of $E$ with respect to this generalized pair is in $(0,1)$.
    In particular, we may apply the proof \autoref{prop:extraction} to this generalized klt pair to find the desired extremal birational model $X'$ extracting exactly $E$. 
\end{proof}

\begin{lemma}\label{lemma:addition}
    Let $k$ be a perfect field of characteristic $p > 5$.
    Let $f \colon X \to Z$ be a Fano type contraction of normal quasi-projective varieties, where $\dim(X)=3$.
    If $Z \to Y$ is a contraction,
    then $\rho(X/Y)=\rho(X/Z)+\rho(Z/Y)$.
\end{lemma}

\begin{proof}
    If $\dim(Z)=0$, the claim is obvious, so we can suppose $\dim(Z) \geq 1$.
    Since $f \colon  X \to Z$ is a Fano type contraction and $p>5$, we have $R^1f_*\mathcal{O}_{X}=0$ by \cite{BK23}*{Theorem 5.1}. Thus we can apply the same proof of \cite{KM92}*{12.1.5.1} (where we replace \cite{KM92}*{12.1.4} with \cite{ABL22}*{Corollary 1.5}) to deduce $\rho(X/Y)=\rho(X/Z)+\rho(Z/Y)$.
\end{proof}

\begin{theorem}\label{thm: fano_model}
    Let $k$ be a perfect field of characteristic $p>5$. 
    Let $(X, \Delta, \mathbf M)$ be a $\qq$-factorial generalized dlt log Calabi--Yau 3-fold.
    Then, there exists a crepant birational model $\phi \colon (X, \Delta, \mathbf M) \dashrightarrow (Y, \Delta_Y, \mathbf M)$ together with a contraction $q \colon Y \to Z$ such that
    \begin{enumerate}
        \item $\Delta_Y^{=1}$ fully supports a $q$-ample $\mathbb{Q}$-divisor;
        \item every log canonical center of $(Y, \Delta_Y, \mathbf M)$ dominates $Z$;
        \item for every $\phi^{-1}$-exceptional divisor $E_Y \subset Y$, we have $E_Y \subset \Delta_Y^{=1}$;
        \item $\phi^{-1}$ is an isomorphism outside $Y \setminus \Delta_Y^{=1}$;
        \item $Z$ is a variety of generalized log Calabi--Yau type.
    \end{enumerate}
    Moreover, $(Y, \Delta_Y, \mathbf M)$ can be chosen to be $\mathbb{Q}$-factorial generalized dlt.
\end{theorem}

\begin{remark}\label{rmk: klt type}
    We observe that property (d) readily implies the weaker property (c), which is stated independently for clarity.
    Properties (c) and (d) in \autoref{thm: fano_model} imply that $\nklt(Y,\Delta_Y,\bM.)=\Delta_Y \sups =1.$ and that $\Delta_Y \geq 0$.
    As $\Delta_Y \sups =1.$ fully supports a $q$-ample $\qq$-divisor, in particular there is an effective divisor $H_Y$ that is $\qq$-Cartier and such that $\Supp(H_Y)=\Supp(\Delta_Y^{=1})=\nklt(Y,\Delta_Y,\bM.)$.
    Then, for $0 < \epsilon \ll 1$, $(Y,\Delta_Y-\epsilon H_Y,\bM.)$ is generalized klt, and it follows from \autoref{lem: pertubation-gen} that $Y$ is of klt type.
\end{remark}

\begin{proof}
For items (a)-(d), we follow the strategy of \cite{FS23}*{Theorem 4.2}.
If $(X,\Delta,\bM.)$ is generalized klt, we take $Y \coloneqq X$ and $Z \coloneqq X$  to conclude.
So, we may assume $\mathrm{Nklt}(X,\Delta,\bM.) \neq \emptyset$.

{\bf Step 0:} In this step, we reduce to the case when the non-klt locus is connected.

Assume $\mathrm{Nklt}(X,\Delta,\bM.)$ is not connected.
Notice that, as $(X,\Delta,\bM.)$ is generalized dlt, this coincides with $\Supp(\Delta \sups =1.)$.
For items (a)-(d), Step 0 in \cite{FS23}*{Theorem 4.2} goes through verbatim, where the existence of the needed MMP is justified by \autoref{lemma:MMP_not_pseff}.
We call $X' \rar Z$ the final model.
By the proof of \cite{FW21}*{Theorem 1.2} and by \autoref{conn_principle}, $(\Delta') \sups =1.$ consists of two disjoint sections for the morphism $X' \rar Z$, and $\bM X'.$ is relatively trivial over $Z$.
By generalized adjunction \cite{FW21}*{Remark 2.5}, we induce a structure of generalized log Calabi--Yau pair on each of the components $T_1$ and $T_2$ of $(\Delta') \sups =1.$. As $T_1 \to Z$ is birational, $Z$ has a generalized log Calabi--Yau structure by \autoref{lem: im-bir-logCY}.
By construction, $X'$ is $\qq$-factorial.
Thus, by (d) and \autoref{lemma:rel_ample}, we can guarantee that $(X',\Delta',\bM.)$ is also generalized dlt.

{\bf Step 1:} We construct a birational map $X \drar X'$ and a contraction $q' \colon X' \rar Z$ satisfying properties (a), (c), and (d) in the statement of the theorem.

This step is identical to Step 1 in \cite{FS23}*{Theorem 4.2}, where the needed $(K_X+\Delta^{<1}+\bM X.)$-MMP exists by \autoref{lemma:MMP_not_pseff} and terminates with a Mori fiber space structure $X' \to Z$.
As $\dim(X')=3$, we have $\dim (Z) \leq 2$.
In particular, if $\dim(Z)=0$, the claim follows.
Therefore, in the rest of the proof, we may assume $\dim(Z) \in \{1,2\}$.
If (b) holds, we conclude. 
Therefore, from now on, we suppose there exists a generalized log canonical center $V$ of $(X',\Delta',\bM.)$ that does not dominate $Z$ and we distinguish the proof according to the dimension of $Z$.

{\bf Step 2:} In this step, we settle properties (a)-(d) in case $\dim(Z)=1$.

Let $V$ be a generalized log canonical center of $(X',\Delta',\bM.)$ that maps onto a closed point $z \in Z$.
If $V$ is a prime divisor, it is contained in $\Supp((\Delta')^{=1,v})$. Since $q' \colon X' \to Z$ is a Mori fiber space onto a curve, we can write $V=(q')^*A$, where $A$ is an effective divisor supported on $\{z\}$.
By construction, $(\Delta')^{=1,h}$ is $q'$-ample.
Since $(\Delta')^{=1,v} \neq 0$, $(\Delta')^{=1,v}$ is the pull-back of an ample $\qq$-divisor on $Z$ and therefore, by taking suitable linear combinations of $(\Delta')^{=1,v}$ and $(\Delta')^{=1,h}$, we conclude that $(\Delta')^{=1}$ fully supports an ample divisor and we replace $Z$ with $\Spec(k)$ to conclude.

Therefore, we may suppose $V$ is not a prime divisor and $(\Delta')^{=1,v}=0$.  
By \autoref{lemma:extraction}, there is an extremal extraction $X'' \to Z$ of a log canonical place $E''$ over $V$.
By construction, $X'$ satisfies (c) and (d).
Then, by \autoref{lemma:rel_ample}, $(\Delta'')^{=1}$ fully supports a $\qq$-divisor $H''$ that is ample over $Z$.
In particular, the model $X''$ satisfies the conditions (a), (c), and (d).
Now, let $D''$ be a general element in $|H''/Z|_{\qq}$.
Then, by condition (d) and Bertini's theorem for ample divisors, $(X'',\Delta''-\epsilon H'' +\epsilon D'',\bM.)$ is generalized klt, where $0 < \epsilon \ll 1$.
Furthermore, this pair is relatively trivial over $Z$.
Then, for $0 < \delta \ll \epsilon$, $(X'',\Delta''-\epsilon H'' +\epsilon D''-\delta E'',\bM.)$ is a generalized klt, since $E''$ is in the support of $\Delta''-\epsilon H''$.
By construction, we have
\begin{equation}\label{eq_perturbe_curve}
    \K X''. + \Delta''-\epsilon H'' +\epsilon D''-\delta E'' + \bM X''.\sim_{\qq,Z} -\delta E''.
\end{equation}
By \autoref{lemma: bz}, this $(-E'')$-MMP over $Z$ terminates with a relatively good minimal model $\widetilde X$.
By construction, we have $\rho(X'/Z)=1$.
Since $f \colon X'' \to X$ is an extremal birational extraction of klt 3-folds and $p>5$, we have $R^1f_*\mathcal{O}_{X''}=0$ by \cite{ABL22}*{Corollary 1.3}.
Thus, we can apply the same proof of \autoref{lemma:addition} to deduce $\rho(X''/Z)=\rho(X''/X')+\rho(X'/Z)=2$.
In particular, the fiber of $X'' \to Z$ consists of 2 irreducible components $E''$ and $F$.
In turn, the $(-E'')$-MMP $X'' \dashrightarrow \widetilde X$ contracts the prime divisor $F$.
In particular, we have $\rho(\widetilde X/Z)=1$.
Thus, $\widetilde \Delta ^{=1}$, which fully supports a divisor that is relatively big over $Z$, fully supports a divisor that is relatively ample over $Z$.
Then, on this new model $\widetilde X$, $\widetilde \Delta ^{=1,v}$ fully supports the pull-back of a divisor that is ample on $Z$.
Hence, we may conclude as in the previous sub-case.
By construction, the final model we obtained is $\qq$-factorial.
Then, by (d) and \autoref{lemma:rel_ample}, we may further assume that the final model of the generalized pair is generalized dlt.

{\bf Step 3:} In this step, we settle properties (a)-(d) in case $\dim(Z)=2$.
We divide this step in various sub-steps.

{\bf Step 3a}: We reduce to the case when $\Supp((\Delta')^{=1,v})$ coincides with the support of the pull-back of an effective $\qq$-Cartier divisor on $Z$.

Assume this is not the case.
By $\rho(X'/Z)=1$, we deduce $(\Delta')^{=1,v}=0$.
By \autoref{lemma:extraction}, we consider an extremal extraction $X'' \to X'$ of a log canonical place $E''$ over $V$ that is vertical over $Z$.
Notice that, by the same argument as in Step 2, $\rho(X''/Z)=\rho(X''/X')+\rho(X'/Z)=2$ holds and that the center of $E''$ in $Z$ is either a curve or a point.
As in Step 2, we may run a $(-E'')$-MMP over $Z$, which terminates with a relatively good minimal model $\widetilde X$.
If the center of $E''$ in $Z$ is a curve $C$, the fiber of $X'' \to Z$ over the generic point of $C$ consists of 2 irreducible components $E''$ and $F$.
In turn, as $-\widetilde E$ is nef over $Z$, the irreducible component $F$ has to be contracted by $X'' \dashrightarrow \widetilde X$.
In particular, we have $\rho(\widetilde X / Z)=1$.
If the center of $E''$ in $Z$ is a point, the relatively ample model resulting from the MMP is a higher birational model $\widetilde Z \to Z$, where the center of $E''$ on $\widetilde Z$ becomes a curve, since $\widetilde E \sim_{\mathbb{Q}, \widetilde Z} 0$. 
By construction, $\rho(\widetilde X /Z) \leq 2$ holds.
In particular, we deduce $\rho(\widetilde X /\widetilde Z)=1$, since $\rho(\widetilde Z / Z) \geq 1$, as $\widetilde Z \to Z$ is not an isomorphism.
Thus, up to replacing the original $X'$ and $Z$ with these models, we may further assume that $(\Delta')^{=1,v}$ is the pull-back of a $\qq$-Cartier divisor on $Z$.
Observe that the indeterminacies of these birational operations are contained in $\Supp((\Delta')^{=1})$, since $X'$ satisfies (d).
In turn, by \autoref{lemma:rel_ample}, (a), (c), and (d) are preserved by this eventual replacements.

{\bf Step 3b:} We reduce to the case when $Z$ has a Mori fiber space structure $Z \to T$. 

First, we endow $Z$ with a generalized log Calabi--Yau structure $(K_Z, B_Z, \mathbf{N})$ and with a generalized pair structure $(Z, \Delta_Z, \mathbf{N})$ with $\Delta_Z<B_Z$. 
Let $Q'$ denote a prime divisor in $\Supp((\Delta')^{=1,v})$, which exists by the reduction in Step 3a.
Note that, since the generic fiber of $X' \to Z$ is a conic, we have that $\Supp((\Delta') \sups =1,h.)$ has either one or two components.
First, assume that $\Supp((\Delta') \sups =1,h.)$ has at least one component that maps birationally to $Z$, and denote it by $S$.
We also write $S^\nu$ for its normalization.
Then, we may perform generalized adjunction for both $(X',\Delta', \bM.)$ and $(X',\Delta' - Q' , \bM.)$ to $S^{\nu}$.
Since both generalized pairs are relatively trivial over $\widetilde Z$, the generalized pair structures obtained on $S^\nu$ can be pushed forward to $Z$ by \autoref{lem: im-bir-logCY}.
Thus, $(Z,B_{Z},\mathbf N)$ and $(Z,\Delta_{Z},\mathbf N)$ are the generalized pairs thus induced on $Z$ from $(X',\Delta', \bM.)$ and $(X',\Delta' - Q',\bM.)$, respectively.
Notice that $B_{Z}=\Delta_{Z} + C$, where the pull-back of $C$ to $X'$ is $Q'$.
Notice that $C$ is effective and $\Supp(C)$ is irreducible.
Now, assume $\Supp((\Delta') ^h)$ is irreducible and maps generically 2 to 1 to $Z$  (which implies $(\Delta') ^ h = (\Delta') \sups =1,h.$).
Set $S \coloneqq \Supp((\Delta') ^h)$, and let $S^\nu$ denote its normalization.
Then, as in the previous case, we may perform generalized adjunction for both $(X',\Delta', \bM.)$ and $(X',\Delta' - Q' , \bM.)$.
In turn, we induce generalized pair structures on $Z$ via \autoref{lem: descent-finite-morph}.
As in the previous case, we have $B_{Z}=\Delta_{Z} + C$.

In both cases, $(Z,B_Z,\mathbf N)$ is generalized log Calabi--Yau, while $K_Z+\Delta_Z+\mathbf N _Z$ is not pseudo-effective.
In turn, we may run a $(K_Z+\Delta_Z+\mathbf N _Z)$-MMP, which ends with a Mori fiber space $Z' \to T$.
Indeed, since $Z$ is a surface and $(Z,\Delta_Z,\mathbf N)$ is generalized log canonical, it follows that $(Z,\Delta_Z)$ is numerically log canonical and in turn log canonical by \cite{Tan18}*{Theorem 4.13}.
Then, by difference, $\mathbf N_Z$ is $\mathbb Q$-Cartier.
Thus, it being the push-forward of a nef divisor on a surface, $\mathbf N_Z$ is itself nef.
Thus, any step of a $(K_Z+\Delta_Z+\mathbf N _Z)$-MMP is also a $(K_Z+\Delta_Z)$-MMP, which can be run by \cite{Tan18}.
As we are in dimension 2, this MMP has finitely many steps and has to terminate.
Since $(K_Z+\Delta_Z+\mathbf N _Z)$ is not pseudo-effective, the MMP indeed terminates with a Mori fiber space.
Notice that this MMP is a $(-C)$-MMP.

Now, we argue that the MMP $Z \to Z'$ can be followed by an MMP on $X'$ such that all the properties of $X'$ and $X' \to Z$ are preserved.
To this end, by finite induction, it suffices to analyze one extremal contraction $Z \to Z_1$.
Let $\Gamma$ be the curve contracted by $Z \to Z_1$ and let $\Xi$ denote its preimage in $X'$.
Since $X' \to Z$ is a Mori fiber space and $X'$ is $\qq$-factorial, $\Xi$ is an irreducible $\qq$-Cartier divisor.
By \autoref{lemma:addition}, we have $\rho(X'/Z_1)=2$.
By assumption, $(\Delta')^{=1}$ fully supports a divisor that is ample over $Z$ and the contraction $Z \to Z_1$ is $C$-positive.
Hence, $Q'$ is the pull-back of a divisor that is relatively ample over $Z_1$.
In turn, we deduce that $(\Delta')^{=1}$ fully supports a divisor that is ample over $Z_1$.
Then, as in Step 2, we can run a $(-Q')$-MMP over $Z_1$, which terminates with a relatively good minimal model $\widetilde X$.
By construction, $\widetilde X$ is $\qq$-factorial, and we have $\rho(\widetilde X /Z_1) \leq 2$.
Furthermore, by $Q' \subset \Supp((\Delta')^{=1})$,
(a), (c), and (d) are preserved.
Since $Z \to Z_1$ contracts $\Gamma$ to a point inside the strict transform of $C$, we have $aQ' + b \Xi\sim_{Z_1} 0$ for some positive integers $a$ and $b$.
In particular, $X' \dashrightarrow \widetilde X$ is a $\Xi$-MMP.
Since the center of $\Xi$ in $Z_1$ is a point, the MMP has to contract $\Xi$ by \cite{Bir12}*{Lemma 3.3}.
Thus, we have $\rho(\widetilde X /Z_1) < \rho(X'/Z_1)=2$.
Thus, we deduce that $\widetilde X \to Z_1$ is a Mori fiber space.

Hence, up to iteratively replacing $X'$ and $Z$ with the models thus constructed, we may further assume that $Z$ is endowed with a Mori fiber space $Z \to T$.
Lastly, we observe that the log Calabi--Yau structure $(Z, B_Z, \mathbf{N})$ on $Z$ is preserved by the surface MMP by push-forward. 

{\bf Step 3c: conclusion.} We show that $X'$ admits a birational transformation to a model $q \colon Y \to Z$ satisfying conditions (a), (b), (c), and (d) of the theorem.

By construction, the divisor $Q'$ is the pull-back from $Z$ of $C$, which is relatively ample over $T$.
Hence, $(\Delta')^{=1}$ fully supports a divisor that is relatively ample over $T$.
By \autoref{lemma:addition}, we deduce  $\rho(X'/T)=2$.
If $T=\Spec(k)$, point (b) holds.
Also, by (d) and \autoref{lemma:rel_ample}, we may further assume that the final model of the generalized pair is generalized dlt.
Similarly, we conclude if $T$ is a curve and all generalized log canonical centers of $(X',\Delta',\bM.)$ dominate $T$.

Thus, we may assume $\dim(T)=1$ and that some generalized log canonical center $V'$ of $(X',\Delta',\bM.)$ is mapped to a point in $T$.
Observe that, since $Z$ dominates $T$, $V'$ does not dominate $Z$.
First assume that $V'$ is a divisor.
Since $X' \to Z$ is a Mori fiber space, $V'$ is the pull-back of a $\qq$-divisor $R$ on $Z$.
In turn, since $Z \to T$ is a Mori fiber space and $V'$ does not dominate $T$, $R$ is a multiple of a fiber of $Z \to T$.
Thus, $V'$ is the multiple of a fiber of $X' \to T$.
Since $T$ is a curve, $V'$ is the pull-back of an ample divisor on $T$, and we may use it to deduce that $(\Delta')^{=1}$ fully supports a globally ample divisor.
Hence, we conclude, and by (d) and \autoref{lemma:rel_ample}, we may further assume that the final model of the generalized pair is generalized dlt.

Thus, we are left with the case when every divisorial component of $(\Delta')^{=1}$ dominates $T$ and $V'$ is a higher codimensional generalized log canonical center.
Let $t$ denote the image of $V'$ in $T$.
By \autoref{lemma:extraction}, we may take an extremal extraction $X'' \to X'$ which extracts a generalized log canonical place $E''$ dominating $V'$.
Notice that this preserves (a), (b), and (d) by \autoref{lemma:rel_ample}.
As in Step 2, we may run a $(-E'')$-MMP relatively over $T$, which terminates with a relatively good minimal model $\widetilde X$.
By \autoref{lemma:addition}, we have $\rho(X''/T)=3$.
Since $- \widetilde E$ is relatively nef and vertical over $T$, it is the multiple of a fiber of $\widetilde X \to T$.
Thus, the MMP $X'' \dashrightarrow \widetilde X$ contracts the strict transform of $X'_t$.
Thus, we have $\rho(\widetilde X /T) \leq 2$.
Since the generic fiber $\widetilde X _\eta$ of $\widetilde X \to T$ coincides with the one of $X' \to T$, we have $\rho(\widetilde X _\eta)=2$.
Since $\widetilde X$ is $\qq$-factorial, we deduce $\rho(\widetilde X /T)=2$.
By construction, $\widetilde X$ coincides with $X'$ over $T \setminus \{t\}$.
In particular, $\widetilde \Delta ^{=1}$ fully supports a divisor that is relatively ample over $T \setminus \{t\}$.
Yet, since we have $\rho(\widetilde X _\eta)=\rho(\widetilde X / T)$, it follows that a divisor that is relatively ample over $T \setminus \{t\}$ is then relatively ample over $T$.
Then, we reduce to the previous sub-case where $\Supp((\Delta')^{=1})$ contains a divisor that is vertical over $T$ and conclude.

{\bf Step 4:} In this step, we verify property (e).

We denote by $q \colon Y \to Z$ the final model constructed in the previous steps and that satisfies properties (a)-(d).
If $\dim(Z)=0$, there is nothing to prove.
If $\dim(Z)=2$, the arguments in Step 3b provide the claim.
Thus, we are left with the case $\dim(Z)=1$.
To this end, let $S'$ denote the normalization of a prime component of $(\Delta') \sups =1,h.$, and we let $(S',B_{ S'},\mathbf L)$ be the generalized log Calabi--Yau pair induced by generalized adjunction.
Let $S' \rar C$ be the Stein factorization of $ S' \rar Z$.
By \autoref{prop:im_gen_CY_Surfaces}, we have $g(C) \leq 1$.
In turn, we have $g(Z) \leq g(C) \leq 1$, and the claim follows.
\end{proof}

We can now prove the existence of rational points on generalized log Calabi--Yau 3-fold pairs with a non-trivial log canonical center.
For the definition of coregularity of a generalized dlt pair, we refer to \cite{FMM}*{Definition 2.4}

\begin{corollary} \label{cor: rat-pts-not-klt}
    Let $p$ be a prime number with $p>5$, and let $q \coloneqq p^e$ for some integer $e \geq 1$.
    Let $(X, \Delta, \bM .)$ be a generalized log Calabi--Yau 3-fold pair, which is not generalized klt, defined over the finite field $\mathbb{F}_q$.
    Let $f \colon X' \to X$ be a generalized dlt modification with $K_{X'}+\Delta'+\bM X'. = f^*(K_X+\Delta+\bM X.)$, whose existence is guaranteed by \autoref{prop: gen_dlt_model}.
    Then, the following statements hold:
    \begin{enumerate}
    	\item if $\coreg(X', \Delta', \bM .') \leq 1$, then   $X(\mathbb{F}_{q}) \neq \emptyset$; and
    	\item if $\coreg(X', \Delta', \bM. ') =2$ and  $q>19$, then $X(\mathbb{F}_{q}) \neq \emptyset$.
    \end{enumerate}
\end{corollary}

\begin{proof}
    By \autoref{cor:rat-pts-finite-fields}, we reduce to prove the statement for the special birational model $Y$ constructed in \autoref{thm: fano_model}, which is of klt type by \autoref{rmk: klt type}.

    As $Z$ is of generalized log Calabi--Yau type by \autoref{thm: fano_model} and it has dimension at most $\coreg(X', \Delta', \bM . ')$, we conclude that $Z(\mathbb{F}_q) \neq \emptyset$ by \autoref{lem: rat_pts_curves} in case (a) (resp. by \autoref{prop:rat_pts_CY?surf}  in case (b)). 
    Let $z \in Z$ be an $\mathbb{F}_q$-rational point.
    If $\dim(Z)=0$, then $\Delta_Y^{=1}$ fully supports an ample $\mathbb{Q}$-divisor and thus we deduce that $Y$ is a 3-fold of Fano type.  By \cite{GNT19}*{Theorem 1.2}, $Y$ has a rational point, concluding.
    Suppose now $\dim(Z) \geq 1$.
    As the fiber $X_{z}$ contains a separably rationally connected geometrically irreducible variety $W$ over $\mathbb{F}_q$ of dimension at most 2 by \autoref{thm: HX-deg-charp}, we apply \autoref{prop: strenghten-Esnault} to conclude.
\end{proof}

\begin{remark}
Note that the bounds on $q$ for the existence of rational points increase with the coregularity. This phenomenon is an incarnation of the more general principle that a pair of coregularity $c$ behaves, under many regards, as a klt pair of dimension $c$.
\end{remark}

\subsection{Generalized klt log Calabi--Yau 3-folds} \label{sect: klt-3-folds}

By \autoref{cor: rat-pts-not-klt}, we may restrict our study to generalized klt log Calabi--Yau 3-folds pairs $(X, \Delta, \bM .)$.

To run an inductive approach, we need to show that the images of 3-folds of generalized log Calabi--Yau type are themselves of generalized log Calabi--Yau type. 
While not able to solve the general case, we have the following partial result, which is sufficient for our purposes.
Our proof crucially relies on the classification of singularities of del Pezzo fibrations and semi-positivity statements in characteristic $p>0$.

\begin{lemma} \label{lemma:extremal-ray}
Let $k$ be an algebraically closed field of characteristic $p > 5$.
Let $(X,\Delta,\bM.)/Z$ be a generalized log canonical 3-fold pair over $k$ such that $X$ is klt.
If $\K X. + \Delta + \bM X.$ is not relatively nef over $Z$, then there exists a $(\K X. + \Delta + \bM X.)$-negative extremal rational curve $\xi$ that is vertical over $Z$.
\end{lemma}

\begin{proof}
    Let $C$ be a curve that is vertical over $Z$ such that $(\K X. + \Delta + \bM X.) \cdot C < 0$.
    Let $A$ be an ample divisor such that $(\K X. + \Delta + \bM X. + A) \cdot C < 0$.
    Then, by \autoref{lem: pertubation-gen}, we may find $\Gamma \sim_{\qq,Z} \Delta + \bM X. + A$ such that $(X,\Gamma)$ is klt.
    Then, the claim follows by \cite{HNT20}*{Theorem 4.6}.
\end{proof}

\begin{proposition} \label{prop: gen_CY_curve_base}
    Let $(X,\Delta, \bM .)$ be a projective $\mathbb Q$-factorial generalized klt log Calabi--Yau 3-fold over a perfect field of characteristic $p>5$.
    Assume that $X$ is geometrically integral, terminal, and admits a $K_X$-Mori fiber space $f \colon X \rightarrow C$ to a curve $C$.
    Then, $C$ has genus at most 1.
\end{proposition}

\begin{proof}
    As $X$ is geometrically integral and the conclusions of the theorem are unaffected by field extensions, we may assume, without loss of generality, that $k$ is algebraically closed.
    Note that $X_{\overline{k}}$ might not be $\mathbb{Q}$-factorial, but $\Delta_{\overline{k}}$ and $\bM X_{\overline{k}}.$ remain $\mathbb{Q}$-Cartier.
    
    Since $(X,\Delta, \bM .)$ is generalized klt, $(X,(1+\epsilon)\Delta, (1+\epsilon) \bM .)$ is generalized klt, where $0 < \epsilon \ll 1$.
    Either $K_X + (1+\epsilon)(\Delta+\bM X.)$ is globally nef, or, by \autoref{lemma:extremal-ray}, there is a rational $(K_X + (1+\epsilon)(\Delta+\bM X.))$-negative extremal curve.
    As $K_X + (1+\epsilon)(\Delta+\bM X.)$ is $f$-nef by assumption, this curve cannot be $f$-vertical.
    Thus, if $K_X + (1+\epsilon)(\Delta+\bM X.)$ is not nef, then $C$ is rational, as it is dominated by a rational curve. 

    Hence, we may assume that $K_X + (1+\epsilon)(\Delta+ \bM X.)$ is globally nef and $f$-ample.
    Equivalently, $\Delta+\bM X.$ is globally nef and $f$-ample.
    By \cite{BT22}*{Theorem 3.7}, the geometric generic fiber $X_{\overline{\eta}}$ is a surface with canonical singularities.
    In particular, it has strongly $F$-regular singularities by \cite{Har98}*{Theorem 1.1}.
    
    Now, let $H$ be an ample divisor on $C$ and fix a rational number $\delta >0$.
    Then, $\Delta+\bM X.+ \delta f^*H$ is ample.
    Hence, by Bertini's theorem, we may choose a general effective divisor $0 \leq \Gamma \sim_{\mathbb Q} \Delta+\bM X.+\delta f^*H$.
    Hence, by Bertini's theorem \cite{SZ13}*{Corollary 6.10}, we may choose a general effective divisor $0 \leq \Gamma \sim_{\mathbb Q} \Delta+\bM X.+\delta f^*H$ such that the geometric generic fiber $(X_{\overline{\eta}}, \Gamma_{\overline{\eta}})$ is strongly $F$-regular.
    Thus, we may apply \cite{Pat14}*{Theorem 3.16} to conclude that $K_{X/C} + \Gamma \sim_{\mathbb Q} K_{X/C}+\Delta+\bM X.+\delta f^*H$ is nef.
    As this argument is valid for any $\delta >0$, we may take the limit as $\delta \to 0$, and we obtain that $K_{X/C}+\Delta+\bM X.$ is nef.
    Then, as $K_X + \Delta + \bM X. \sim_{\mathbb Q} 0$, we obtain that $-f^*K_C$ is nef, and hence we may conclude that the genus of $C$ is at most 1.
\end{proof}

If the image of a 3-fold $X$ of generalized log Calabi--Yau type is a normal surface $Z$, in general, we are not able to show that $Z$ is of generalized log Calabi--Yau type.
Nevertheless, we can show strong geometric constraints on $Z$ up to universal homeomorphism.

\begin{proposition}\label{prop: point_base_conic_bundle_gen}
    Let $k$ be a perfect field of characteristic $p>5$.
    Let $(X,\Delta, \bM X.)$ be a geometrically integral $\mathbb{Q}$-factorial generalized klt 3-fold log Calabi--Yau pair over $k$.
    Let $f \colon X \to S$ be a $K_X$-Mori fiber space of relative dimension 1.
    Then $S$ has  $\mathbb{Q}$-factorial rational singularities and there exists a purely inseparable finite morphism $\psi \colon T \to S$ such that $K_T$ is either numerically trivial or it is not pseudo-effective.
\end{proposition}

\begin{proof}
    By \cite{BK23}*{Theorem 27}, $S$ has rational singularities and thus is $\mathbb{Q}$-factorial by \cite{kk-singbook}*{Proposition 10.9}.
    If $K_S$ is numerically trivial or not pseudo-effective, we conclude.
    So from now, we assume that $K_S$ is pseudo-effective and not numerically trivial.
    We will show that there exists a purely inseparable cover $\psi \colon T \rar S$ such that $-K_T$ is pseudo-effective.

    Since $\bM X.$ is pseudo-effective and $f$ is a Mori fiber space, $\bM X.$ is either $f$-ample or $\bM X. \sim_{\qq} f^*L$ for some pseudo-effective $\qq$-Cartier $\qq$-divisor $L$ on $S$ by \cite{HNT20}*{Theorem 5.4}.

    {\bf Case 1:} In this case, we assume $-(K_X+\Delta) \sim_\qq \bM X. \sim_\qq f^*L$ for some pseudo-effective $\qq$-Cartier $\qq$-divisor $L$ on $S$.
  
    The generic fiber is a regular conic by \cite{kk-singbook}*{Lemma 10.6} and, as $p>2$, it is smooth by \cite{BT22}*{Lemma 2.17}. 
    By \cite{Wit21}*{Theorem 3.4}, we may find a purely inseparable finite morphism $\psi \colon T \rar S$ such that
    \begin{equation}\label{eq:jakub_epsilon_easy}
        -\psi^* L \sim_{\qq} t K_{T} + (1-t)\psi^*K_S + E,
    \end{equation}
    where $t \in [0,1]$ and $E \geq 0$.
    If $t=0$, \autoref{eq:jakub_epsilon_easy} implies $\psi^* K_S \sim_\qq -\psi^*L -E$.
    In particular, $-K_S$ is pseudo-effective, thus leading to a contradiction.
    Then, $t>0$ holds.
    In turn, \autoref{eq:jakub_epsilon_easy} implies $t K_T \sim_\qq -\psi^*L-E-(1-t)\psi^*K_S$.
    In particular, $-K_T$ is pseudo-effective, and the claim follows.

    {\bf Case 2:} In this case, we assume that $\bM X.$ is $f$-ample.
    
    Let $\pi \colon X' \rar X$ be a log resolution where $\bM.$ descends.
    Also, fix an ample $\qq$-divisor $A$ on $S$.
    Since $f$ has relative dimension 1, $\pi$ is an isomorphism over the generic point of $S$.
    Thus, $\bM X'.$ is globally nef and $f$-big.
    Then, for every integer $n > 0$, $\bM X'. + \frac 1 n g^* A$ is nef and big, where we set $g \coloneqq f \circ \pi$.
    Write $\K X'. + \Delta' + \bM X'. = \pi^*(\K X. + \Delta + \bM X.)$.
    In particular, $(X',\Delta')$ is sub-klt.
    Then, by the same arguments as in the proof of \autoref{lem: pertubation-gen}, we may find $0 \leq \Gamma'_n \sim_{\qq} \bM X'. + \frac 1 n g^* A$ such that $(X',\Delta'+\Gamma'_n)$ is sub-klt.
    In particular, since $\bM X'. + \frac 1 n g^* A$ is nef and big, there is an effective $\qq$-divisor $F'_n$ such that, for every positive integer $l$, we may find an ample divisor $H'_l$ such that $\bM X'. + \frac 1 n g^* A \sim_\qq H'_l+ \frac{1}{l}F'_n$.
    In particular, the support of $\frac{1}{l}F'_n$ is independent of $l$.
    Thus, we may choose $l$ large enough such that the following holds:
    If the geometric generic fiber $(X'_{\overline{\eta}},\Delta'_{\overline{\eta}})$ is sub-klt, then so is $(X'_{\overline{\eta}},\Delta'_{\overline{\eta}}+\frac{1}{l}F'_{n,\overline{\eta}})$.
    If we fix a regular closed fiber $X'_{z}$ over a closed point $z \in S$, by \cite{7authors}*{Theorem 2.17}, we can choose $H'_l$ such that its support intersects $X'_{z}$ transversally.
    Since the field $k(z)$ is a finite extension of $k$, it is perfect, so it follows that the intersection between the support of $H'_l$ and $X_z$ is smooth.
    In particular, the support of $H'_l$ is separable over $S$.
    Furthermore, by ampleness, we may assume that its support does not share components with the supports of $\Delta'$ and $F'_n$.
    In particular, if $(X'_{\overline{\eta}},\Delta'_{\overline{\eta}}+\frac{1}{l}F'_{n,\overline{\eta}})$ is sub-log canonical, then so is $(X'_{\overline{\eta}},\Delta'_{\overline{\eta}}+\frac{1}{l}F'_{n,\overline{\eta}}+H'_{l,\overline{\eta}})$.
    Then, by defining $\Gamma'_n \coloneqq H'_l+\frac{1}{l}F'_n$ for $l$ large enough and $H'_l$ general in its $\qq$-linear series, and setting $\Gamma_n \coloneqq \pi_* \Gamma'_n$, we get that $(X,\Delta+\Gamma_n)$ is klt.
    Notice that $\K X. + \Delta + \Gamma_n \sim_{\qq} {\frac 1 n} f^* A$.

    The generic fiber is a regular conic by \cite{kk-singbook}*{10.3} and, as $p>2$, it is smooth by \cite{BT22}*{Lemma 2.17}. 
    By \cite{Wit21}*{Theorem 3.4}, we may find a purely inseparable finite morphism $\psi_n \colon T_n \rar S$ such that
    \begin{equation}\label{eq:jakub_epsilon}
        \frac{1}{n} \psi_n^* A \sim_{\qq} t_n K_{T_n} + (1-t_n)\psi_n^*K_S + E_n,
    \end{equation}
    where $t_n \in [0,1]$ and $E_n \geq 0$.

    \begin{claim}\label{claim:epsilon}
        Up to considering a subsequence, the morphism $\psi_n \colon T_n \rar S$ does not depend on $n$, $\{E_n\}_{n \in \mathbb N}$ is convergent up to $\qq$-linear equivalence, and $\{ t_n \}_{n \in \mathbb N}$ is a converging sequence with limit $t_{\infty} \in [0,1]$.
    \end{claim}

    We will provide the proof of \autoref{claim:epsilon} after this current proof.
    By passing to a subsequence, by \autoref{claim:epsilon}, we may write $\psi_n=\psi$ and $T_n= T$.
    Furthermore, as $E_n$ depends continuously on $\frac 1 n$, we may set $E \coloneqq \lim_{n \to \infty} E_n$.
    While the limit $E$ may not be well defined as a divisor, it is well defined as a numerical equivalence class, and as such it will be understood in the following.
    In particular, $E$ is pseudo-effective.
    First, assume that $t_{\infty}=0$ in \autoref{eq:jakub_epsilon} (in case $T=S$, we declare $t_{\infty}=0$ and treat it in this case).
    Then, we have ${\frac 1 n} \psi^*A \sim_{\qq} \psi^*K_S + E_n$.
    Since $E_n$ depends continuously on ${\frac 1 n}$, we may take the limit as $n \to \infty$, and we obtain $K_S \equiv -E$.
    Since $K_S$ is pseudo-effective and not numerically trivial and $E$ is pseudo-effective, we reach a contradiction.
    Therefore, $t_{\infty}>0$ holds true.
    Again by a limiting argument, we obtain $t_{\infty} K_T \equiv - (1-t_{\infty})\psi^*K_S - E$, thus concluding that $-K_T$ is pseudo-effective.
\end{proof}

\begin{proof}[{Proof of \autoref{claim:epsilon}}]
    For the setup and notation of this proof, we refer to the proof of \autoref{prop: point_base_conic_bundle_gen} where \autoref{claim:epsilon} is invoked.
    Let $\eta$ denote the generic point of $S$ and $\overline{\eta}$ its geometric generic point.

    First, assume that, for infinitely many $n$, $(X_{\overline{\eta}},\Delta_{\overline{\eta}}+\Gamma_{n,\overline{\eta}})$ is log canonical.
    Then, we may apply \cite{Wit21}*{Proposition 3.2}.
    In particular, we have $T_n=S$ and $\psi_n= \mathrm{id}_S$ for infinitely many $n$.
    Furthermore, for such $n$, we have $E_n \sim_{\qq} -K_S + \frac{1}{n} A$.

    Now, assume $(X_{\overline{\eta}},\Delta_{\overline{\eta}}+\Gamma_{n,\overline{\eta}})$ is not log canonical for infinitely many $n$.
    Thus, we have that $(X'_{\overline{\eta}},\Delta'_{\overline{\eta}}+\frac{1}{l} F'_{n,\overline{\eta}})$ is not log canonical by the above description of $H'_{l}$.
    We claim there exists, up to passing to a subsequence, a prime divisor $\Sigma$ in $\Supp(\Delta)$ such that some of the coefficients of the irreducible components of $\Sigma_{\overline \eta}$ appear in $(\Delta'_{\overline{\eta}}+\frac{1}{l} F'_{n,\overline{\eta}})$ with coefficient strictly larger than 1. 
    First, we suppose $(X_{\overline{\eta}},\Delta_{\overline{\eta}})$ is log canonical.
    Then $F_{n}$, which denotes the push-forward of $F'_n$ to $X$, contains an irreducible component $G_n$ of the support $\Delta^{=1, h}$, and we choose $\Sigma=G_n$. 
    Notice that $\Supp(\Delta^{=1, h})$ is independent of $n$, so up to a subsequence, we may assume that  $\Sigma=G_n$ is independent of $n$.
    If instead $(X_{\overline{\eta}},\Delta_{\overline{\eta}})$ is not log canonical, for every $n$ we choose as $\Sigma$ the component of $\Supp(\Delta^{h})$ that forces $(X_{\overline{\eta}},\Delta_{\overline{\eta}})$ not to be log canonical and, up to passing to a subsequence, we can suppose that $\Sigma$ is indipendent of $n$.
 
   For more details, see \cite{Wit21}*{proof of Theorem 3.4} (in the notation of \emph{op. cit.}, $\Sigma$ corresponds to $S$).
    In particular, by \cite{Wit21}*{proof of Theorem 3.4}, $T_n$ and $\psi_n$ are independent of $n$ for infinitely many $n$, as the base change is uniquely determined by the horizontal divisor $\Sigma$.

    To conclude, we need to study $E_n$ and $t_n$.
    Let $a_0$ denote the coefficient of $\Sigma$ in $\Delta$.
    Then, by the proof of \autoref{lem: pertubation-gen}, we have that the coefficient of $\Sigma$ in $\Delta+\Gamma_n$ is $a_n=a_0+o(n^{-1})$ (this may be achieved by taking $\sigma$ sufficiently small for $n$ fixed in the proof of \autoref{lem: pertubation-gen}).
    In the proof of \cite{Wit21}*{Theorem 3.4}, the value of $t_n$ is determined by the linear equation $a_n=\frac{1}{p^d}t_n+(1-t_n)$, where $p$ is the characteristic of the field, and $d$ is an integer determined uniquely by $\Sigma$ ($S$ in \emph{op.~cit.}).
    Thus, as we have $a_n=a_0+o(n^{-1})$, the claim about $\{ t_n\}_{n \in \mathbb N}$ immediately follows.
    By the proof of \cite{Wit21}*{Theorem 3.4}, $E_n$ is determined, up to $\qq$-linear equivalence, by the adjunction to the normalization of $\Sigma$ of the divisors $t_n(K_X+ \frac{1}{p^d}\Sigma)$, $(1-t_n)(K_X + \Sigma)$, and $(\Delta+\Gamma_n-a_n\Sigma)$.
    Since $\{ t_n\}_{n \in \mathbb N}$ is convergent and the divisors $K_X+ \frac{1}{p^d}\Sigma$ and $K_X + \Sigma$ do not depend on $n$, the first two contributions to $E_n$ are convergent (even as $\qq$-divisors).
    Then, by construction, we have
    \[
    \Delta+\Gamma_n-a_n\Sigma \sim_{\qq} \Delta-(a_0-o(n^{-1}))\Sigma + \bM X. + \frac{1}{n} f^*A,
    \]
    and the claim follows.
\end{proof}

\begin{corollary} \label{cor: model-purely-insep}
    Let $k$ be a perfect field of characteristic $p > 5$.
    Let $(X,\Delta, \bM X.)$ be a geometrically integral $\mathbb{Q}$-factorial generalized klt 3-fold log Calabi--Yau pair over $k$.
    Let $f \colon X \to S$ be a $K_X$-Mori fiber space of relative dimension 1.
    Then $S$ has rational $\mathbb{Q}$-factorial singularities and there exists a purely inseparable finite morphism $\psi \colon T \to S$ such that one of the following holds:
    \begin{enumerate}
        \item $T$ is birational to a smooth surface with numerically trivial canonical class;
        \item $T$ is a geometrically rational surface; or
        \item $T$ admits a contraction $p \colon T \to C$ onto a curve of genus $g(C) = 1$ such that $-K_T$ is big over $C$.
    \end{enumerate}
\end{corollary}

\begin{proof}
    Let $\psi \colon T \to S$ be the inseparable cover constructed in \autoref{prop: point_base_conic_bundle_gen}.
    As $S$ has rational singularities by \autoref{prop: point_base_conic_bundle_gen},
    then $T$ has $W \mathcal{O}$-rational singularities by \cite{PZ21}*{Proposition 3.11}. We distinguish the two cases of \autoref{prop: point_base_conic_bundle_gen} according to the positivity of $K_T$. 
    
    Suppose $K_T \sim_{\mathbb{Q}} 0$.
    If $T$ has canonical singularities, then we are done by taking the minimal resolution.
    If $T$ has worse than canonical singularities, we consider a log resolution $\pi \colon T' \to T$.
    As $T$ has $W\mathcal{O}$-rational singularities, the exceptional locus of $\pi$ is a tree of  rational curves by \cite{NT20}*{Proposition 2.23}.
    As $S$ has $\mathbb{Q}$-factorial singularities by \cite{kk-singbook}*{Proposition 10.9},
    $T$ is also $\mathbb{Q}$-factorial by \cite{Tan18}*{Lemma 2.5} and thus the same proof of \cite{BS23}*{Lemma 2.15} shows that $T$  is a geometrically rational surface.
    
    Suppose now $K_T$ is not pseudo-effective.
    If $T$ is geometrically rational, we conclude.
    If not, let $T' \to T$ be the minimal resolution.
    As $K_{T'}$ is not pseudo-effective, by \cite{Tan18} we can run a $K_{T'}$-MMP, which gives a birational contraction $T' \to T''$ of smooth projective surfaces where $T''$ admits a Mori fiber space structure $T'' \to C$. As $T'$ is not geometrically rational, by the classification of surfaces \cite{Mum69}, the base $C$ of the Mori fiber space is a curve of genus $g(C) \geq 1$. In particular, $T'$ admits a fibration $T' \to C$ with rationally chain connected fibers onto a curve of genus $g(C) \geq 1$. To conclude the proof of (c), we are thus left to show $g(C) = 1$.
    As $T$ has $W \mathcal{O}$-rational singularities, by the rigidity lemma \cite{Deb01}*{Proposition 1.14}, the fibration $T' \to C$ descends to a fibration $T \to C$ as all exceptional divisors of $T' \to T$ are rational curves by \cite{NT20}*{Proposition 2.23}.
    As $T \to S$ is a universal homeomorphism, we have an isomorphism $\Pic(S)_{\mathbb{Q}} \to \Pic(T)_{\mathbb{Q}}$ by \cite{Kee99}*{Lemma 1.4.(3)}.
    Applying \cite{Kee99}*{Lemma 1.4.(2)}, there exists a commutative diagram 
    $$
    \xymatrix{
    T \ar[r]^{\psi} \ar^{}[d] & S \ar[d]^{}   \\
    C \ar[r]^{\varphi} & D
    }
    $$
    where $S \to D$ is a contraction.
    Consider the composite morphism $X \to S \to D$.
    By \autoref{prop: gen_CY_curve_base}, the curve $D$ has genus at most 1.
    Since $\phi$ is a universal homeomorphism, also $C$ has genus at most 1, thus concluding.
\end{proof}
    
Now, we prove the existence of rational points on 3-folds of generalized log Calabi--Yau type whose underlying variety has negative Kodaira dimension.

\begin{theorem} \label{prop: strictly-klt-3-CY}
Let $p$ be a prime number with $p>5$, and let $q \coloneqq p^e$ for some integer $e \geq 1$.
Let $(X, \Delta, \bM X.)$ be a geometrically integral generalized log Calabi--Yau 3-fold over $\mathbb{F}_q$.
Suppose that $\kappa(X)=-\infty$. 
If we have $q>19$, then $X(\mathbb{F}_{q}) \neq \emptyset$ holds.
\end{theorem}

\begin{proof}
    By passing to a terminalization, we can suppose that $X$ is terminal and $\mathbb{Q}$-factorial.
    Then, we can run a $K_X$-MMP, which will end with a Mori fiber space $X' \to B$ by the assumption that $\kappa(X)=-\infty$, the non-vanishing theorem \cite{XZ}, and the fact that $X$ is terminal.
    By \autoref{cor:rat-pts-finite-fields}, we may replace $X$ with $X'$.
    Thus, without loss of generality, we may assume that $X$ is terminal, $\mathbb Q$-factorial, and it admits a Mori fiber space $f \colon X \rightarrow B$.
    Now, we proceed by cases depending on $\dim (B)$.

    If $\dim (B)=0$, then $X$ is a 3-fold of Fano type and hence it has a rational point by \cite{GNT19}*{Theorem 1.2}.
    
    If $\dim (B)=1$, then we have $g(B) \leq 1$ by \autoref{prop: gen_CY_curve_base}.
    In particular, $B$ has a rational point $b$ by \autoref{lem: rat_pts_curves}. 
    By \autoref{thm: HX-deg-charp}, the fiber $X_b$ contains a separably rationally connected subvariety that is geometrically irreducible over $k(b)$. 
    Then, we conclude $X_b$ has a rational point by \autoref{prop: strenghten-Esnault}.
    
    If $\dim (B) =2$, consider the purely inseparable finite morphism $T \to S$ given by \autoref{cor: model-purely-insep}.
    Let (a), (b), and (c) denote the cases in \autoref{cor: model-purely-insep}.
    In cases (a) and (b), $T$ has a rational point by \autoref{cor: surf_bir_rat_pts}.
    In case (c), we proceed as in the proof of \autoref{prop:rat_pts_CY?surf}.
    In particular, $T$ fibers over a curve of genus at most 1, which has a rational point by \autoref{lem: rat_pts_curves}.
    In turn, $T$ has a rational point by combining \cite{BT22}*{Proposition 2.18} and the Chevalley--Warning's theorem \cite{Ser73}*{\S~I.I.2 Theorem 3} on an appropriate smooth model of $T$.
    Then, as $T$ admits a rational point and maps to $S$, it follows that also $S$ has a rational point.
    Again by \autoref{thm: HX-deg-charp} and \autoref{prop: strenghten-Esnault} we conclude the fiber over an $\mathbb{F}_{q}$-rational point on $S$ has a $\mathbb{F}_q$-rational point, concluding the proof.
\end{proof}

We conclude this section by showing the Albanese morphism is surjective.

\begin{proposition} \label{prop: gen_klt_CY_albanese}
    Let $k$ be a perfect field of characteristic $p>5$.
    Let $(X, \Delta, \bM X.)$ be a generalized log Calabi--Yau 3-fold pair defined over $k$.
    Suppose that $\kappa(X)=-\infty$. 
    Then, the Albanese morphism of $X$ is surjective.
\end{proposition}

\begin{proof}
    By \autoref{lem: im_alb}, we can follow the same proof of \autoref{prop: strictly-klt-3-CY} and we can suppose that $X$ is terminal and it admits a Mori fiber space $f \colon X \to B$. 
    As the fibers of $f$ are rationally chain connected, by the rigidity lemma \cite{Deb01}*{Proposition 1.14}, it suffices to show that the Albanese morphism of $B$ is surjective.

    If $\dim (B)=0$, this is immediate.
    If $\dim (B)=1$, the geometric genus $g(B)$ is at most 1 by \autoref{prop: gen_CY_curve_base}, and we conclude.

    If $\dim (B) =2$, let $T \to B$ be the purely inseparable finite morphism given by \autoref{cor: model-purely-insep}.
    As $T \to S$ is a universal homeomorphism, we have a factorization of a power of the Frobenius morphism $F^{(e)} \colon S \to T \to S$ and thus we have that $\alb_S$ is surjective if and only if $\alb_T$ is surjective by \autoref{lem: im_alb}.
    A case-by-case analysis easily implies that $\alb_T$ is surjective.
    By \autoref{lem: im_alb}, we may replace $T$ with the minimal model of its minimal resolution.
    Let (a), (b), and (c) denote the cases in \autoref{cor: model-purely-insep}.
    In cases (a) and (b), $\alb_T$ is surjective by \autoref{prop: alb_surf}.
    In case (c), $T$ fibers over a curve $C$ with $g(C) \leq 1$ and $\alb_T$ coincides with $\alb_C$.
    Since $g(C) \leq 1$, we conclude.
\end{proof}

\subsection{Canonical Calabi--Yau 3-folds with non-trivial Albanese}

In virtue of the following immediate consequence of \autoref{prop: strictly-klt-3-CY}, we are left to study the case of
canonical $K$-trivial 3-folds.

\begin{lemma}\label{remark:kod_dim}
    Let $k$ be a perfect field of characteristic $p>5$, and let $(X,\Delta, \bM.)$ be a generalized log Calabi--Yau 3-fold over $k$.
    If $\kappa(X) \neq -\infty$ holds, then $X$ is canonical and $\K X. \sim_\qq 0$ holds. 
\end{lemma}
\begin{proof}    
By \autoref{prop:extraction}, we may find a birational morphism $\phi \colon X' \rar X$ such that $(X',\Delta',\bM.)$ is a generalized dlt modification of $(X,\Delta, \bM.)$ and $X'$ is terminal and $\qq$-factorial.
    By construction, we have $\kappa(\K X'.)=\kappa(X')=\kappa(X) \geq 0$.
    Since $\K X'. + \Delta' + \bM X'. \sim_{\qq} 0$ and $\kappa(\K X'.) \geq 0$ hold, it follows that $\Delta'=0$ and $\bM X'.=0$ hold.
    In turn, we have $\Delta=0$ and $\bM. =0$, where the latter is understood as an identity of b-divisors, and the claim follows.
\end{proof}

Although we are not able to solve the general case concerning $K$-trivial 3-folds, we construct rational points in the case that the Albanese morphism is not trivial.

We start by showing that the Albanese morphism of canonical $K$-trivial 3-folds is surjective, building on the recent work of Ejiri and Patakfalvi \cite{EP23}.

\begin{proposition} \label{prop: surj_albanese_K_trivial}
	Let $k$ be a perfect field of characteristic $p>0$.
    Let $X$ be a canonical 3-fold with $K_X \equiv 0$. 
    Then
    $\alb_X \colon X \to \Alb_X$ is surjective.
\end{proposition}

\begin{proof}
	We can assume $k$ to be algebraically closed.
    Let $\phi \colon X' \to X$ be a terminalization.
    By \cite{GNT19}*{Theorem 4.8}, the $\phi$-exceptional locus is rationally chain connected.
    Thus, by the rigidity lemma \cite{Deb01}*{Proposition 1.14}, we have $\alb_{X'}=\alb_X \circ \phi$.
    Thus, we may assume that $X$ is $\mathbb{Q}$-factorial and terminal.
    If the image of $\alb_X$ is a point, then by the universal property of the Albanese morphism, $\alb_X$ is surjective.
    If the image of $\alb_X$ has positive dimension, the generic fiber of $\alb_X$ is regular, as terminal 3-folds have isolated singularities \cite{kk-singbook}*{Theorem 2.29}.
    Thus, we conclude by \autoref{thm: surj_EP}.
\end{proof}

Now, we show the surjectivity of the Albanese morphism for generalized log Calabi--Yau 3-fold pairs.

\begin{theorem}\label{thm: surj_albanese}
    Let $k$ be a perfect field of characteristic $p > 5$.
    Let $(X, \Delta, \bM .)$ be a generalized log Calabi--Yau 3-fold over $k$.
    Then $\alb_X$ is surjective.
\end{theorem}

\begin{proof}
    First, assume that $(X,\Delta,\bM.)$ is generalized klt.
    If $\kappa(X)= -\infty$ holds, the claim follows from \autoref{prop: gen_klt_CY_albanese}.
    Otherwise, we apply \autoref{prop: surj_albanese_K_trivial}.
    Therefore, in the remainder of the proof, we may assume that $(X,\Delta,\bM.)$ is not generalized klt.
    
    By \autoref{lem: im_alb}, we may replace $(X,\Delta,\bM.)$ with a dlt modification.
    Thus, we may assume that $(X,\Delta,\bM.)$ is $\qq$-factorial generalized dlt.
    Now, let $(Y,\Delta_Y,\bM.)$ be a crepant birational model as in \autoref{thm: fano_model}, and let $q \colon Y \rar Z$ the corresponding fibration to a lower dimensional variety of generalized log Calabi--Yau type.
    First, we show that $\alb_Y$ is surjective, and then we argue that this implies that so is $\alb_X$.

    By \autoref{prop: gen_CY_curve_base} and \autoref{prop: alb_surf}, $\alb_Z$ is surjective.
    By \autoref{thm: fano_model}, there exists an effective $q$-ample divisor $H_Y$ such that $\Supp(H_Y)=\Supp(\Delta_Y^{=1})$.
    Furthermore, by \autoref{rmk: klt type}, $(Y,\Delta_Y-\epsilon H_Y,\bM.)$ is generalized klt for $0 < \epsilon \ll 1$.
    Fix such $\epsilon$.
    Let $H_Z$ be an ample divisor on $Z$ such that $H_Y+q^*H_Z$ is ample.
    Then, we apply \autoref{lem: pertubation-gen} to the generalized klt pair $(Y,\Delta_Y-\epsilon H_Y,\bM.)$ and the ample divisor $\frac{\epsilon}{2}(H_Y+q^*H_Z)$.
    Thus, we may find an effective divisor $\Gamma_Y$ such that $(Y,\Gamma_Y)$ is klt and $-(K_Y+\Gamma_Y)$ is $q$-ample.
    Then, by \cite{GNT19}*{Theorem 4.1}, every fiber of $q$ is rationally chain connected.
    In turn, by \cite{Deb01}*{Proposition 1.14}, $\alb_Y=\alb_Z \circ q$ holds.
    In particular, as $\alb_Z$ is surjective, so is $\alb_Y$.
    
    Thus, we are left with comparing $\alb_X$ and $\alb_Y$.
    By construction, $X$ is a klt variety.
    Furthermore, by \autoref{rmk: klt type}, $Y$ is of klt type.
    Consider a common resolution $W$.
    Then, by \cite{GNT19}*{Theorem 4.8}, the exceptional loci of $\alpha \colon W \rar X$ and $\beta \colon W \rar Y$ are rationally chain connected.
    Then, by \cite{Deb01}*{Proposition 1.14}, $\alb_X \circ \alpha=\alb_W=\alb_{Y} \circ \beta$ holds and the claim follows.
\end{proof}

To show the existence of a rational point for a $K$-trivial 3-fold with non-trivial Albanese, we need to understand the geometric generic fiber of the Albanese morphism.
We start by showing the geometric rationality of geometrically non-normal regular $K$-trivial surfaces over an imperfect field of characteristic $p\geq 5$. 

\begin{proposition} \label{prop: K-trivial-not-norm-rat}
Let $k$ be a field of characteristic $p>0$ and let $X$ be a geometrically integral regular projective surface over $k$ with $K_X \equiv 0$.
If $p > 3$ and $X$ is not geometrically normal, then $X$ is a geometrically rational surface.
\end{proposition}

\begin{proof}
We can suppose $k$ to be separably closed.
Consider $f \colon Y \coloneqq (X \times_k \overline{k})^\nu \to X$, where $(X \times_k \overline{k})^\nu \rar X \times_k \overline{k}$ is the normalization.
By \cite{PW22}*{Theorem 1.1}, there exists an effective non-trivial Weil divisor $C$ on $Y$ such that $K_Y + (p-1)C \sim f^* K_X$.
If $\pi \colon Z \to Y$ is the minimal resolution, then $K_Z$ is not pseudo-effective.
Thus we have that either
\begin{enumerate}
\item[(i)] $Z$ is rational; or
\item[(ii)] there exists a fibration $g \colon Z \to B$ to a curve $B$ of genus $g(B)>0$ whose generic fiber is rational by \cite{Mum69}.
\end{enumerate} 

Now, we will exclude case (ii).
As $k$ is separably closed, $X_{\overline{k}} \to X$ is a universal homeomorphism and so is the composition $Y \to X$ by \cite{Tan18}*{Lemma 2.2}.
Then, we conclude that $Y$ is $\mathbb{Q}$-factorial by \cite{Tan18}*{Lemma 2.5}.
As $k$ is imperfect (and thus not contained in the algebraic closure of a finite field), by \cite{Tan14}*{Theorem 3.20}, we conclude that all $\pi$-exceptional curves are $g$-vertical.
Write $K_Z+\Gamma=\pi^*K_Y$ and $\pi^*C=\pi_*^{-1}C+E$, where $E$ and $\Gamma$ are effective $\mathbb{Q}$-Cartier $\mathbb{Q}$-divisors.
Then, for a general fiber $F$ of $g$, we have
\begin{equation}
\begin{split}    
0 &= \pi^* f^* K_X \cdot F\\
&= \pi^*(K_Y+(p-1)C) \cdot F\\
&= K_Z \cdot F + \Gamma \cdot F+ (p-1)((\pi_*^{-1} C + E) \cdot F) \\
&=-2+(p-1)(\pi_*^{-1} C \cdot F),
\end{split}
\end{equation}
where the first equality follows from the fact that $K_X$ is numerically trivial, and the last equality follows from the fact that $F$ is a conic and that $\Gamma$ is vertical over $B$ and hence its intersection with a general fiber is 0.
From the equation and the smoothness of $Z$, it follows that $\pi_*^{-1} C \cdot F$ is a positive integer.
Then, we deduce that $p \leq 3$, reaching a contradiction.
\end{proof}

\begin{remark}
The assumption on the characteristic in \autoref{prop: K-trivial-not-norm-rat} is sharp.
For $p \in \{2,3\}$, let $E$ be the regular non-smooth planar cubic curve defined by $\left\{y^2z=x^3-tz^3 \right\}$ over $\mathbb{F}_p(t)$.
If $C$ is a smooth elliptic curve over $\mathbb{F}_p(t)$ then $X \coloneqq E \times_{\mathbb{F}_p(t)} C$ is a regular surface with $K_{X} \sim 0$ and the normalized base change $(X \times_{k} \overline{k})^{\nu}$ is $\mathbb{P}^1 \times C$, which is not rational.
\end{remark}

\begin{proposition}\label{prop:last_statement}
    Let $p$ be a prime number with $p > 3$, and let $q \coloneqq p^e$ for some integer $e \geq 1$.
    Let $X$ be a geometrically integral canonical 
    projective $K$-trivial 3-fold over $\mathbb{F}_q$ such that its Albanese variety has positive dimension.
    If we have $q>19$, then $X(\mathbb{F}_{q}) \neq \emptyset$ holds.
\end{proposition}

\begin{proof}
    As in the proof of \autoref{prop: surj_albanese_K_trivial}, we can suppose $X$ to be terminal and $\qq$-factorial.
	By \autoref{prop: surj_albanese_K_trivial}, we have $\dim(\Alb_X) \leq 3$. 
	Let $X \to Z \to \Alb_X$ be the Stein factorization of the Albanese morphism $\alb_X$.
    By
    \autoref{thm: surj_EP},
    the finite morphism $Z \to \Alb_X$ is purely inseparable and thus there exists $e>0$ such that there is a factorization of the Frobenius morphism $F^e \colon \Alb_X \to Z \to \Alb_X$. In this case, by \cite{Lan55} we have $\Alb_X(\mathbb{F}_q) \neq \emptyset$ and thus we deduce $Z(\mathbb{F}_q) \neq \emptyset$. 
    
    If $\dim (Z)=3$, then $X \to Z$ is birational, thus it is separable and we conclude by 
    \autoref{thm: surj_EP}
    that $Z \simeq \Alb_X$ holds and $X \to \Alb_X$ is birational. As $K_X$ is trivial, and in particular nef over $\Alb_X$, $\Alb_X$ is smooth, we conclude $X = \Alb_X$ by the negativity lemma.
    Therefore we conclude $X(\mathbb{F}_q) \neq \emptyset$ by \autoref{cor:rat-pts-finite-fields}.
    
    Suppose $\dim (\Alb_X)=2$.
    As $p>3$, the geometric generic fiber of $X \to Z$ is a smooth elliptic curve by \cite{PW22}*{Theorem 1.4}.
    In particular, the geometric generic fiber is integral, so $X \to Z$ is a separable contraction.
    By \autoref{thm: surj_EP} (or \cite{CWZ23}*{Theorem 8.1}) and \autoref{thm: iso_EP},
    we conclude that $Z=\Alb_X$ and all fibers are geometrically integral smooth curves of genus 1.
    Let $z \in Z$ be an $\mathbb{F}_q$-rational point and let $F$ be the fiber over $z$.
    As $X$ is Cohen--Macaulay by \cite{BK23}, we can apply \cite{PZ}*{Theorem 9.1} to conclude that $\alb_X$ is isotrivial.
    Therefore, $F_z$ is geometrically an elliptic curve, and we conclude by \cite{Lan55} that $X$ has an $\mathbb{F}_q$-rational point.

    Finally, suppose $\dim (\Alb_X) =1$. Then $Z=\Alb_X$ holds and $X \to Z$ is isotrivial 
    by \autoref{thm: iso_EP}
    and the geometric generic fiber is reduced by \cite{Sch10}. 
    Let $z$ be an $\mathbb{F}_q$-rational point on $Z$ and let $F$ be the closed fiber over $z$, which is geometrically integral by the isotriviality of $\alb_X$.
    First, suppose that the geometric generic fiber is normal.
    By the isotriviality of $\alb_X$, then $F$ is normal and $K_F \sim_{\mathbb{Q}} 0$.
    If $F$ is log canonical, we conclude $F(\mathbb{F}_{q}) \neq \emptyset$ by \autoref{prop:rat_pts_CY?surf}. 
    If $F$ is not log canonical, then by \cite{FW21}*{Theorem 1.2}, the non-klt locus is geometrically connected.
    In particular, the non-klt locus defines an $\mathbb{F}_q$-rational point.
    If the geometric generic fiber is not normal, then $F$ is geometrically rational by \autoref{prop: K-trivial-not-norm-rat}.
    By the isotriviality of $\alb_X$, $F$ is an integral non-normal geometrically rational surface.
    By \cite{Kol96}*{IV.3.2.5}, $F$ is rationally connected and thus $F(\mathbb{F}_q) \neq \emptyset$ holds by \autoref{prop: strenghten-Esnault}.
\end{proof}

In the proofs of \autoref{prop: surj_albanese_K_trivial} and \autoref{prop:last_statement}, we need some technical generalizations of \cite{EP23}*{Theorem 1.3} and \cite{EP23}*{Theorem 1.6}.
These generalizations are reported here below for completeness.
We do not claim any authorship for these statements, as they appear as courtesy of Ejiri and Patakfalvi.
The statements below will then appear in an update to \cite{EP23}.

\begin{theorem}[{cf. \cite{EP23}*{Theorem 1.3}}] \label{thm: surj_EP}
    Let $(X,\Delta)$ be a projective pair with $-(K_X+\Delta)$ nef. 
    Let $X \to \Alb_X$ be the Albanese morphism and $X \rar Z \rar \Alb_X$ its Stein factorization.
    Suppose that $K_X+\Delta$ is a $\mathbb{Z}_{(p)}$-Weil divisor and that the generic fiber $(X_\eta, \Delta_\eta)$ is strongly $F$-regular and $p$ does not divide the Cartier index of $K_{X_\eta}+\Delta_\eta$.
    Then the following hold:
    \begin{enumerate}
        \item the morphism $\alb_X$ is surjective;
        \item the morphism $Z \to \Alb_X$ is purely inseparable; and
        \item if $X \to Z$ is separable, then $Z=\Alb_X$ holds.
    \end{enumerate}
\end{theorem}

\begin{theorem}[{cf. \cite{EP23}*{Theorem 1.6}}] \label{thm: iso_EP}
    Let $(X,\Delta)$ be a projective pair with $-(K_X+\Delta)$ nef. 
    Let $\alb_X\colon X \to \Alb_X$ be the Albanese morphism and assume that $\alb_X(X)$ has dimension 1.
    Suppose that $K_X+\Delta$ is a $\mathbb{Z}_{(p)}$-Weil divisor and that the generic fiber $(X_\eta, \Delta_\eta)$ is strongly $F$-regular and $p$ does not divide the Cartier index of $K_{X_\eta}+\Delta_\eta$.
    If $-(K_X+\Delta)$ is semi-ample, then $\alb_X$ is an isotrivial algebraic fiber space.
\end{theorem}

\bibliographystyle{amsalpha}
\bibliography{refs}
	
\end{document}